\documentclass{amsproc}

\usepackage{graphs}
  \usepackage{color}

\newtheorem{theorem}{Theorem}[section]

\theoremstyle{definition}
\newtheorem{definition}[theorem]{Definition}
\newtheorem{example}[theorem]{Example}

\newtheorem{proposition}[theorem]{Proposition}
\newtheorem{corollary}[theorem]{Corollary}
\theoremstyle{remark}
\newtheorem{remark}[theorem]{Remark}

\numberwithin{equation}{section}

\def\rank{{\rm rank}}

\def\wt{\widetilde}
\def\ol{\overline}
\def\ov{\overline}

\def\wh{\widehat}

\newcommand{\N}{{\mathbb N}}
\newcommand{\Z}{{\mathbb Z}}

\def\ov{\overline}

\newcommand{\om}{\omega}

\def\Rk{{\mathcal R}}

\def\qed{\hfill$\square$}

\def\sig{\sigma}

\numberwithin{equation}{section}
\begin{document}

\title{Bratteli diagrams: structure, measures, dynamics}

\author{S. Bezuglyi}
\address{Department of Mathematics, Institute for Low Temperature Physics, Kharkiv 61103, Ukraine}
\curraddr{Department of Mathematics, University of Iowa, Iowa City,
52242 IA, USA}
\email{bezuglyi@gmail.com}

\author{O. Karpel}
\address{Department of Mathematics, Institute for Low Temperature Physics, Kharkiv 61103, Ukraine}
\email{helen.karpel@gmail.com}

\date{}

\dedicatory{To the memory of Ola Bratteli}


\begin{abstract} This paper is a survey on general (simple and non-simple) Bratteli diagrams which focuses on the following topics: finite and infinite tail invariant measures on the path space $X_B$ of a Bratteli diagram $B$, existence of continuous dynamics on $X_B$ compatible with tail equivalence relation, subdiagrams and measure supports. We also discuss the structure of Bratteli diagrams,  orbit equivalence and full groups, homeomorphic measures.
\end{abstract}

\maketitle
\tableofcontents

\section{Introduction}\label{introduction}
This paper is devoted to {\em Bratteli diagrams}, the object that is widely used for constructions  of  transformation models in various dynamics.
It is difficult to overestimate the significance of Bratteli diagrams for the study of dynamical systems. A class of graduated infinite graphs, later called Bratteli diagrams, was originally introduced by Bratteli \cite{Br72} in his breakthrough article on the classification of approximately finite (AF) $C^*$-algebras.

It turned out that the close ideas developed by Vershik in the study of sequences of measurable partitions led to a realization of any ergodic automorphism of a standard measure as a transformation acting on a path space of a graph (afterwards called a Vershik map) \cite{vershik:1981}, \cite{vershik:1982}.  The dynamical systems obtained in this way are called  Bratteli-Vershik dynamical systems. During the last two decades, Bratteli diagrams turned out to be a very powerful and productive tool for the study of dynamical systems not only on a measure space but also on Cantor and Borel spaces \cite{herman_putnam_skau:1992}, \cite{BDK06}. By a Cantor dynamical system we mean a pair $(X,T)$ consisting of a Cantor set $X$ and a homeomorphism $T : X \to X$. The results proved in \cite{herman_putnam_skau:1992}  build a bridge between Cantor dynamics and Bratteli diagrams. It was proved that any minimal Cantor dynamical system $(X, T)$  is realized as a Bratteli-Vershik homeomorphism defined on the path space $X_B$ of a Bratteli diagram $B$. The diagrams arising in this way have a nice property: they are {\em simple}.

Our goal is to show that a large part of results proved in the context of Cantor minimal dynamical systems remains true for a much wider class of aperiodic homeomorphisms of a Cantor set. First of all, every aperiodic homeomorphism admits its realization as a Vershik map on a {\em non-simple} Bratteli diagram \cite{medynets:2006}. Moreover, the non-simple stationary and finite rank  Bratteli diagrams correspond to substitution dynamical systems and expansive homeomorphisms (the same as in minimal dynamics) (see \cite{durand_host_skau}, \cite{BKM09}, \cite{DM}). On the other hand, the problem of the existence of a continuous dynamics on a non-simple Bratteli diagram is much harder. It turns out that there are Bratteli diagrams that cannot serve as Bratteli-Vershik dynamical systems. But nevertheless one can consider the {\em tail equivalence relation} on such diagrams. This relation determines a kind of dynamical system on the path space.  Our discussions of this issue are mostly based on  \cite{medynets:2006}, \cite{BKY14}, \cite{BY}, \cite{JQY14}. Also there are no general results about classification of aperiodic homeomorphisms with respect to the orbit equivalence relation. In contrast to minimal case, where a number of nice theorems were proved (see \cite{giordano_putnam_skau:1995}, \cite{glasner_weiss:1995-1},  \cite{GPS99},
\cite{giordano_putnam_skau:2004}, \cite{GMPS08}, \cite{GMPS10}, \cite{HKY12}), we are aware of only some sort of negative results which show that the invariants used in minimal case do not work, in general, for aperiodic homeomorphisms \cite{S.B.O.K.} (the only exclusion is full groups \cite{Med}).

The main reason why Bratteli diagrams are convenient to use for the study of homeomorphisms $T: X\to X$  is the fact that  various properties of $T$  become more transparent  when one deals with corresponding Bratteli-Vershik dynamical systems. This observation is related to $T$-invariant measures and their supports, to minimal components of $T$, structure of $T$-orbit, etc.
 In particular, the problem of finding all ergodic $T$-invariant measures (the extreme points of the Choquet simplex $M(X,T)$) and their supports for a given $(X, T)$ is traditionally a central one in the theory of dynamical systems, especially for specific interesting examples of  homeomorphisms $T$. But being considered in general settings, this problem looks rather vague, and there are very few universal results that can be applied to  a given homeomorphism $T$. But for Bratteli-Vershik realization of $T$, we are going to discuss some natural methods for the study of $M(X,T)$ based on the structure of the underlying diagram. Moreover, these methods work even for Bratteli diagrams that do not support any Vershik map. To emphasize the difference between simple and non-simple Bratteli diagrams, we remark that for an aperiodic homeomorphism $T$ the simplex  $M(X,T)$ may contain infinite measures. This is impossible for minimal dynamical systems.

There are important applications of Bratteli diagrams. One of them is the theory of countable dimension groups. For every such a group $G$ there exists a Bratteli diagram whose $K^0$ group is order isomorphic to $G$. In case of simple Bratteli diagrams this group can be realized as the quotient group of $C(X, \Z)$ by $T$-boundaries. Then $T$-invariant measures are in one-to-one correspondence with traces of the dimension group. Dimension groups  have a large repertoire of results and examples, and sometimes it is convenient to translate the original problems of measures to this context. For instance, we consider the notion of good measures (defined by Akin \cite{Akin2}) for stationary Bratteli diagrams. It is interesting to see how this notion can be extended to traces of dimension groups  \cite{BH}. Handelman studied this notion and other amazing properties of traces in his recent works. Another application is considered in \cite{BJ} where some representations of the Cuntz-Krieger algebras are constructed by stationary Bratteli diagrams.

We should say that the list of papers devoted to different aspects of Bratteli diagrams is very long. It is impossible to cover (or even mention) most of them. We have to restrict our choice of the material for this work to several topics that are clearly mentioned in the title. Unfortunately, we were unable to discuss relations between Bratteli diagrams and operator algebras otherwise it would doubled the size of the paper. Also the notion of dimension groups and impressive results establishing links to Bratteli diagrams are not included in this survey.

There are several recent surveys that are focused on other important directions of the study of Bratteli diagrams (see, for example, \cite{durand:2010}, \cite{P10}, \cite{S00}).  Our survey has minor intersections with them.

\medskip

\textbf{Notation}
\\

\begin{tabular}{ll}

$(X,T)$ &  Cantor dynamical system \\
$B = (V, E)$ &  Bratteli diagram \\
$B = (V, E, \omega)$ &  ordered Bratteli diagram \\
$X_B$ & path space of a Bratteli diagram\\
$\varphi_\om, \varphi$ &  Vershik map on ordered diagram $B$\\
$\mathcal{R}$ &  tail equivalence relation\\
$Orb_T(x)$ & $T$-orbit of $x$\\
\end{tabular}

\section{Fundamentals of Cantor dynamics and Bratteli diagrams}\label{Preliminaries}

This section contains the basic definitions and facts about Cantor dynamical systems and Bratteli diagrams. They are the main objects of our consideration in this paper.

\subsection{Cantor dynamics} By definition, a {\em Cantor set (space)} $X$ is  a zero-dimensional compact metric space without isolated points. The topology on $X$ is generated by a countable family of clopen subsets. It is called the clopen topology. All such Cantor sets are homeomorphic.

A {\em homeomorphism} $T : X \to X$ is a continuous bijection. Denote
$Orb_T(x) := \{T^n (x) \; |\; n \in \mathbb{Z}\}$; the set $Orb_T(x)$ is called the {\em orbit} of $x\in X$ under action of $T$ (or simply $T$-orbit). We consider here only {\em aperiodic} homeomorphisms $T$, i.e., for every $x$ the set $Orb_T(x)$  is countably infinite. In fact, some definitions and facts make sense for arbitrary homeomorphisms. We focus on the case of aperiodic homeomorphisms for convenience mostly.

A homeomorphism $T : X \to X$ is called {\em minimal} if for every $x \in X$ the set $Orb_T(x)$ is dense. A {\em minimal Cantor system}  is a pair $(X,T)$ where $X$ is a Cantor space, and $T \colon X \rightarrow X$ is a minimal homeomorphism. Any (aperiodic) homeomorphism $T$ of a Cantor set has a {\em minimal component} $Y$: this is a $T$-invariant closed non-empty subset $Y$ of $X$ such that $T|_Y$ is minimal on $Y$.

Given a minimal Cantor system $(X,T)$ and a clopen  $A\subset X$, the first return function  $r_A(x) = \min\{n \geq 1 : T^n(x) \in A\}$ is a well defined  continuous integer-valued function with domain $A$. Then $T_A(x) = T^{r_A}(x)$ is a homeomorphism of $ A$, and $(A, T_A)$ is  called the \textit{induced} Cantor minimal system.

There are several notions of equivalence for  Cantor dynamical systems. We give these definitions for single homeomorphisms of Cantor sets.

\begin{definition}
Let $(X,T )$ and $(Y, S)$ be two aperiodic Cantor systems. Then

(1) $(X,T)$ and $(Y,S)$ are \textit{conjugate} (or \textit{isomorphic}) if there exists a homeomorphism $h \colon X \rightarrow Y$ such that $h \circ T = S \circ h$.

(2) $(X,T)$ and $(Y,S)$ are \textit{orbit equivalent} if there exists a homeomorphism $h \colon X \rightarrow Y$ such that $h(Orb_T(x)) = Orb_S(h(x))$ for every $x \in X$. In other words, there exist functions $n, m \colon X \rightarrow \mathbb{Z}$ such that for all $x \in X$, $h \circ T(x) = S^{n(x)} \circ h(x)$ and $h \circ T^{m(x)} = S \circ h(x)$. The functions $n,m$ are called the orbit cocycles associated to $h$.

(3)  $(X,T)$ and $(Y,S)$ are \textit{strong orbit equivalent} if they are orbit equivalent and each of the corresponding orbit cocycles $n, m$ has at most one point of discontinuity.

(4) $(X,T)$ and $(Y,S)$ are \textit{Kakutani equivalent} if they both have clopen subsets such that the corresponding induced systems are conjugate.

(5) $(X,T)$ and $(Y,S)$ are \textit{Kakutani orbit equivalent} if they both have clopen subsets such that the corresponding induced systems are orbit equivalent.
\end{definition}

Given a Cantor dynamical system $(X,T)$, a Borel measure $\mu$ on $X$ is called $T$-\textit{invariant} if $\mu(TA) = \mu(A)$ for any Borel set $A$. Let $M(X,T)$ be the set of all invariant measures. It is well known that  $M(X,T)$ is a Choquet simplex whose extreme points are $T$-ergodic measures. This simplex includes probability measures (when $\mu(X) =1$)  and infinite measures (when $\mu(X) = \infty$). We observe that infinite invariant measures may exist for aperiodic homeomorphisms; in minimal dynamics this is impossible. If $M(X,T) =\{\mu\}$, then $T$ is called {\em uniquely ergodic}.

It is not hard to see that every Borel measure $\mu$ on a Cantor set is completely determined by its values on clopen sets. This means that there is a one-to-one correspondence between $\mu$ and the collection of numbers $S(\mu) = \{\mu(A) : A\ \mbox{clopen}\} \subset [0,1]$.

It follows from \cite{herman_putnam_skau:1992}, \cite{giordano_putnam_skau:1995}, and \cite{medynets:2006} that any minimal (and even aperiodic) Cantor dynamical system $(X,T)$ admits a realization as a Bratteli-Vershik dynamical system $(X_B, \varphi_B)$ acting on a path space $X_B$ of a Bratteli diagram (see Section~\ref{BV}). Thus, the study of $T$-invariant measures is reduced to the case of measures defined on the path space of a Bratteli diagram. The advantage of this approach is based on the facts that (i) any such a measure is completely determined by its values on cylinder sets of $X_B$, and (ii) there are simple and explicit formulas for measures of cylinder sets. Especially transparent this method works for stationary and finite rank Bratteli diagrams, simple and non-simple ones \cite{BKMS_2010}, \cite{BKMS_2013}.

It is worth pointing out that the study of measures on a Bratteli diagram is a more general problem than that in  Cantor dynamics. This observation follows from the existence of Bratteli diagrams that do not support any continuous dynamics on their path spaces which is compatible with the tail equivalence relation. The first example of such a Bratteli diagram was given  in \cite{medynets:2006}; a more comprehensive coverage of this subject can be found in \cite{BKY14} and \cite{BY} (we discuss this stuff below in Section~\ref{VM}). If a Bratteli diagram does not admit a Bratteli-Vershik homeomorphism, then we have to work with the {\em tail equivalence relation} $\mathcal R$ on $X_B$ and study measures invariant with respect to $\mathcal R$.

\subsection{Bratteli diagrams}
 A {\it Bratteli diagram} is an infinite graph $B=(V,E)$ such that the vertex
set $V =\bigcup_{i\geq 0}V_i$ and the edge set $E=\bigcup_{i\geq 0}E_i$
are partitioned into disjoint subsets $V_i$ and $E_i$ where

(i) $V_0=\{v_0\}$ is a single point;

(ii) $V_i$ and $E_i$ are finite sets, $\forall i \geq 0$;

(iii) there exist $r : V \to E$ (range map $r$) and $s : V \to E$ (source map $s$), both from $E$ to $V$, such that $r(E_i)= V_{i+1}$, $s(E_i)= V_{i}$, and
$s^{-1}(v)\neq\emptyset$, $r^{-1}(v')\neq\emptyset$ for all $v\in V$
and $v'\in V\setminus V_0$.

The set of vertices $V_i$ is called the $i$-th level of the
diagram $B$. A finite or infinite sequence of edges $(e_i : e_i\in E_i)$
such that $r(e_{i})=s(e_{i+1})$ is called a {\it finite} or {\it infinite
path}, respectively. For $m<n$, $v\, \in V_{m}$ and $w\,\in V_{n}$, let $E(v,w)$ denote the set of all paths $\overline{e} = (e_{1},\ldots, e_{p})$
with $s(\ol e) = s(e_{1})=v$ and $r(\ol e) = r(e_{p})=w$. If $m < n$ let $E(m,n)$ denote all paths whose source belongs to $V_m$ and whose range belongs to $V_n$.
For a Bratteli diagram $B$,
let $X_B$ be the set of infinite paths starting at the top vertex $v_0$.
 We endow $X_B$ with the topology generated by cylinder sets
$[\ol e]$ where $\ol e= (e_0, ... , e_n)$,  $n \in \mathbb N$, and
$[\ol e]:=\{x\in X_B : x_i=e_i,\; i = 0, \ldots, n\}$.
 With this topology, $X_B$ is a 0-dimensional compact metric space.
By assumption, we will consider only such Bratteli diagrams $B$ for which $X_B$ is a {\em Cantor set}, that is $X_B$ has no isolated points.  Letting $|A|$ denote the cardinality of the set $A$, this means that for every $(x_0,x_1,\ldots)\in X_B$
and every $n\geq 1$ there exists $m>n$ such that $|s^{-1}(r(x_m))|>1$.

Given a Bratteli diagram $B$, the $n$-th {\em incidence matrix}
$F_{n}=(f^{(n)}_{v,w}),\ n\geq 0,$ is a $|V_{n+1}|\times |V_n|$
matrix such that $f^{(n)}_{v,w} = |\{e\in E_{n+1} : r(e) = v, s(e) = w\}|$
for  $v\in V_{n+1}$ and $w\in V_{n}$.  Every vertex $v \in V$ is connected with $v_0$ by a finite path, and the set of $E(v_0, v) $ of all such paths is finite. If $h^{(n)}_v = |E(v_0, v)|$, then for all $n \geq 1$
\begin{equation}
h_v^{(n+1)}=\sum_{w\in V_{n}}f_{v,w}^{(n)}h^{(n)}_w \ \mbox{or} \ \
h^{(n+1)}=F_{n}h^{(n)}
\end{equation}
where $h^{(n)}=(h_w^{(n)})_{w\in V_n}$. The numbers $h_w^{(n)}$ are usually called heights (see Section~\ref{BV}).

We define the following important classes of Bratteli diagrams that we work with in this article.

\begin{definition}\label{rank_d_definition} Let $B$ be a Bratteli diagram.
\begin{enumerate}
\item
We say that $B$ has \textit{finite rank} if  for some $k$, $|V_n| \leq k$ for all $n\geq 1$.
\item
Let $B$ have finite rank. We say that $B$ has \textit{rank $d$} if  $d$ is the smallest integer such that $|V_n|=d$ infinitely often.
\item
We say that $B$ is {\em simple} if  for any level
$n$ there is $m>n$ such that $E(v,w) \neq \emptyset$ for all $v\in
V_n$ and $w\in V_m$. Otherwise, $B$ is called {\em non-simple}.
\item
 We say that $B$ is  \textit{stationary} if $F_n = F_1$  for all $n\geq 2$.
\end{enumerate}
\end{definition}

Let $x =(x_n)$ and $y =(y_n)$  be two paths from $X_B$. It is said that $x$ and $y$ are {\em tail equivalent} (in symbols,  $(x,y) \in \mathcal R$)  if there exists some $n$ such that $x_i = y_i$ for all $i\geq n$. Since $X_B$ has no isolated points, the $\mathcal R$-orbit of any point $x\in X_B$ is infinitely countable. The diagrams with infinite $\mathcal R$-orbits are called {\em aperiodic}. Note that a Bratteli diagram is simple if the tail equivalence relation $\mathcal R$ is minimal.

In order to illustrate the above definitions, we give an example of a Bratteli diagram.

\unitlength = 0.5cm
 \begin{center}
 \begin{graph}(17,15)
\graphnodesize{0.4}
 \roundnode{V0}(8.5,14)
\roundnode{V11}(0.0,9) \roundnode{V12}(4.25,9)
\roundnode{V13}(8.5,9) \roundnode{V14}(12.75,9)
\roundnode{V15}(17.0,9)
\roundnode{V21}(0.0,5) \roundnode{V22}(4.25,5)
\roundnode{V23}(8.5,5) \roundnode{V24}(12.75,5)
\roundnode{V25}(17.0,5)
\roundnode{V31}(0.0,1) \roundnode{V32}(4.25,1)
\roundnode{V33}(8.5,1) \roundnode{V34}(12.75,1)
\roundnode{V35}(17.0,1)
\edge{V11}{V0} \edge{V12}{V0} \edge{V13}{V0} \edge{V14}{V0}
\edge{V15}{V0}
\edge{V21}{V11}\bow{V21}{V11}{0.06} \edge{V21}{V12} \edge{V22}{V11} \edge{V22}{V12}
\edge{V23}{V13} \edge{V23}{V14} \edge{V24}{V13} \edge{V24}{V14}
\edge{V25}{V11}\bow{V25}{V11}{0.06} \edge{V25}{V13} \bow{V25}{V15}{0.06}
\bow{V25}{V15}{-0.06} \edge{V25}{V15}
\edge{V31}{V21} \edge{V31}{V22} \bow{V31}{V22}{0.06}\edge{V32}{V21} \edge{V32}{V22}
\edge{V33}{V23}\bow{V33}{V23}{0.06} \edge{V33}{V24} \edge{V34}{V23} \edge{V34}{V24}
\edge{V35}{V21} \edge{V35}{V23} \bow{V35}{V25}{0.06}
\bow{V35}{V25}{-0.06} \edge{V35}{V25}

\end{graph}
\centerline{.\ .\ .\ .\ .\ .\ .\ .\ .\ .\ .\ .\ .\ .\ .\ .\ }
\\
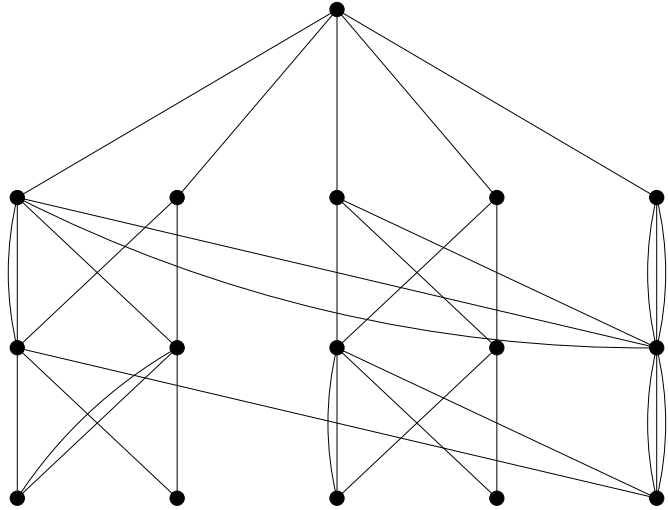
Fig. 1. Example of a Bratteli diagram

 \end{center}
\medskip

This diagram is a non-simple  finite rank Bratteli diagram that has exactly two minimal components (they are clearly seen).

We will  constantly use the \textit{telescoping} procedure for a Bratteli diagram:

\begin{definition} \label{telescoping_definition}
Let $B$ be a Bratteli diagram, and $n_0 = 0 <n_1<n_2 < \ldots$ be a strictly increasing sequence of integers. The {\em telescoping of $B$ to $(n_k)$} is the Bratteli diagram $B'$, whose $k$-level vertex set $V_k'$ is $V_{n_k}$ and whose incidence matrices $(F_k')$ are defined by
\[F_k'= F_{n_{k+1}-1} \circ \ldots \circ F_{n_k},\]
where $(F_n)$ are the incidence matrices for $B$.
\end{definition}

Roughly speaking, in order to telescope a Bratteli diagram, one takes a subsequence of levels $\{n_k\}$ and considers the set $E(n_k, n_{k+1})$  of all finite paths between the  levels $\{n_k\}$ and $\{n_{k+1}\}$ as edges of the new diagram.
In particular, a Bratteli diagram $B$ has rank $d$ if and only if there is
a telescoping $B'$ of $B$ such that $B'$ has exactly $d$ vertices
at each level. When telescoping diagrams, we often do not specify to which levels $(n_k)$ we telescope, because it suffices to know that such a sequence of levels exists.

Two Bratteli diagrams are {\em isomorphic} if they are isomorphic as graded graphs.
Define an equivalence relation $\sim$ on the set of Bratteli diagrams generated by isomorphism and telescoping. One can show that $B_1 \sim B_2$ if there exists a Bratteli diagram $B$ such that  telescoping of $B$ to odd levels yields, say, $B_1$ and telescoping to even levels yields $B_2$.

In order to avoid consideration of some trivial cases, we will assume that the following \textit{convention} always holds: {\em our Bratteli diagrams are not unions of two or more disjoint subdiagrams.}

\subsection{Ordered Bratteli diagrams and Vershik maps} The concept of an ordered Bratteli diagram is crucial for the existence of dynamics on the path space of a Bratteli diagram.

\begin{definition}\label{order_definition} A Bratteli diagram $B=(V,E) $ is called {\it ordered} if a linear order `$>$' is defined on every set  $r^{-1}(v)$, $v\in
\bigcup_{n\ge 1} V_n$. We use  $\om$ to denote the corresponding partial
order on $E$ and write $(B,\om)$ when we consider $B$ with the ordering $\om$. Denote by $\mathcal O_{B}$ the set of all orderings on $B$.
\end{definition}

Every $\omega \in \mathcal O_{B}$ defines the  \textit{lexicographic}
ordering on the set $E(k,l)$ of finite paths between
vertices of levels $V_k$ and $V_l$:  $(e_{k+1},...,e_l) > (f_{k+1},...,f_l)$
if and only if there is $i$ with $k+1\le i\le l$, such that $e_j=f_j$ for $i<j\le l$
and $e_i> f_i$.
It follows that, given $\om \in \mathcal O_{B}$, any two paths from $E(v_0, v)$
are comparable with respect to the lexicographic ordering generated by $\om$.  If two infinite paths are tail equivalent, and agree from the vertex $v$ onwards, then we can compare them by comparing their initial segments in $E(v_0,v)$. Thus $\om$ defines a partial order on $X_B$, where two infinite paths are comparable if and only if they are tail equivalent.

   \begin{definition}
We call a finite or infinite path $e=(e_i)$ \textit{ maximal (minimal)} if every
$e_i$ is maximal (minimal) amongst the edges
from $r^{-1}(r(e_i))$.
\end{definition}

Notice that, for $v\in V_i,\ i\ge 1$, the
minimal and maximal (finite) paths in $E(v_0,v)$ are unique. Denote
by $X_{\max}(\om)$ and $X_{\min}(\om)$ the sets of all maximal and
minimal infinite paths in $X_B$, respectively. It is not hard to see that
$X_{\max}(\om)$ and $X_{\min}(\om)$ are \textit{non-empty closed subsets} of
$X_B$; in general, $X_{\max}(\om)$ and $X_{\min}(\om)$ may have
interior points. For a finite rank Bratteli diagram $B$, the sets $X_{\max}(\om)$
and $X_{\min}(\om)$ are always finite for any $\om$, and if $B$ has rank $d$,
then each of them have at most $d$ elements
(\cite{BKM09}). For an  aperiodic Bratteli diagram $B$, we see that   $X_{\max}(\om)\cap X_{\min}(\om) =\emptyset$.

We say that an  ordered Bratteli diagram $(B, \om)$ is  \textit{properly ordered} if the sets $X_{\max}(\om)$ and $X_{\min}(\om)$ are singletons.

If we denote by $\mathcal O_B(j)$ the set of all orders on $B$ which have $j$ maximal and $j$ minimal paths, then, in this notation,  $\mathcal O_B(1)$ is the set of proper orders on $B$.

Let $(B,\omega)$ be an ordered Bratteli diagram, and suppose that $B'=(V',E')$ is the telescoping of $B$ to levels $(n_k)$. Let $v' \in V'$ and suppose that the two edges $e_1',$ $e_2'$, both with range $v'$, correspond to the finite paths $e_1$, $e_2$ in $B$, both with range $v$. Define the order $\om'$ on $B'$ by $e_1'<e_2'$ if and only if $e_1<e_2$. Then $\om' $ is called the {\em lexicographic order generated by $\om$} and is denoted by $\om'=L(\om)$.

 It is not hard to see that if $\om' = L(\om)$, then
$$
|X_{\max}(\om)| = |X_{\max}(\om')|,\ \ |X_{\min}(\om)| = |X_{\min}(\om')|.
$$

A Bratteli diagram $B $ is called \textit{ regular} if for any ordering
$\omega \in \mathcal O_{B}$ the sets $X_{\max}(\omega)$ and
$X_{\min}(\omega)$ have empty interior.

In particular, finite rank Bratteli diagrams are automatically regular, and if all incidence matrix entries of $B$ are at least 2, then $B$ is regular. We consider here only regular Bratteli diagrams.

We will need the notion of the language associated to an ordered Bratteli diagram. If $V$ is a finite alphabet, let $V^{+}$ denote the set of nonempty
words over $V.$ We use the notation $W' \subseteq W$ to indicate that $W'$ is a subword of $W$.
 If $W_1, W_2, \ldots ,W_n$, are words, then we let  $\prod_{i=1}^{n}W_i$ refer to their concatenation.

Let $\om$ be an order on a Bratteli diagram $B$.  Fix a vertex $v\in
V_n$ and some level $m<n$, consider the set $E(V_m,v) = \bigcup_{v'\in V_{m}}E(v',v)$ of all finite paths between vertices of level $m$ and $v$. This set can be ordered by $\om$: $E(V_m, v) = \{e_{1},\ldots e_{p}\}$ where $e_{i}<e_{i+1}$ for
$1\leq i\leq p-1$. Define the word $w(v,m,n):=s(e_{1})s(e_{2})\ldots
s(e_{p})$ over the alphabet $V_{m}$. If $W=v_{1}\ldots v_{r} \in
V_{n}^{+}$, let $w(W,n-1, n):= \prod_{i=1}^{r}w(v_{i},n-1,n).$

\begin{definition}\label{language_definition}
The {\em level-$n$ language $\mathcal
L(B,\om, n)$ of $(B,\om)$ }  is
$$\mathcal L(B,\om,n):=
\{W :  \,\, W\subseteq w(v,n,N), \,\, \mbox{ for
some } v\in V_{N},  N>n\}\, .
$$
\end{definition}

If $B$ has strict rank $d$, then each of the level-$n$ languages can be defined on a common alphabet $V$, and in this case we have the {\em language} of $B$
 $$ \mathcal L(B, \om):=\limsup_{n}\mathcal L(B,\om, n)\, .
$$

The idea to use an order on a Bratteli diagram to define a transformation acting
on the path space $X_B$ was firstly developed  by Vershik \cite{vershik:1981}, and then it was applied in many papers (see, e.g. \cite{durand:2010}, \cite{giordano_putnam_skau:1995},
\cite{herman_putnam_skau:1992})

\begin{definition}\label{VershikMap}
Let $(B, \omega)$ be an ordered, regular Bratteli
diagram. We say that $\varphi = \varphi_\omega : X_B\rightarrow X_B$
is a {\it (continuous) Vershik map} if it satisfies the following conditions:

(i) $\varphi$ is a homeomorphism of the Cantor set $X_B$;

(ii) $\varphi(X_{\max}(\omega))=X_{\min}(\omega)$;

(iii) if an infinite path $x=(x_0,x_1,\ldots)$ is not in $X_{\max}(\omega)$,
then $\varphi(x_0,x_1,\ldots)=(x_0^0,\ldots,x_{k-1}^0,\overline
{x_k},x_{k+1},x_{k+2},\ldots)$, where $k=\min\{n\geq 1 : x_n\mbox{
is not maximal}\}$, $\overline{x_k}$ is the successor of $x_k$ in
$r^{-1}(r(x_k))$, and $(x_0^0,\ldots,x_{k-1}^0)$ is the minimal path
in $E(v_0,s(\overline{x_k}))$.
\end{definition}

If $\om$ is an ordering on $B$, then one can always define the map
$\varphi_0$ that maps $X_B \setminus X_{\max}(\om)$ onto $X_B
\setminus X_{\min}(\om)$ according to (iii) of Definition
\ref{VershikMap}. The question about the existence of the Vershik map is
equivalent to that of an extension of  $\varphi_0 : X_B \setminus
X_{\max}(\om) \to X_B \setminus X_{\min}(\om)$ to a homeomorphism
of the entire set $X_B$.  If $\om$ is a proper ordering, then
$\varphi_\om$ is a homeomorphism. In particular any simple Bratteli diagram has a Vershik map. For a finite rank Bratteli diagram $B$, the
situation is simpler than for a general Bratteli diagram because the sets
$X_{\max}(\om)$ and $X_{\min}(\om)$ are finite.

\begin{definition}\label{Good_and_bad}  Let $B$ be a Bratteli diagram
$B$.  We say that an ordering $\om\in \mathcal O_{B}$ is \textit{perfect}
if $\om$ admits a Vershik map $\varphi_{\om}$ on $X_B$.
Denote by  $\mathcal P_B $ the set of
all perfect orderings  on $B$.
\end{definition}

We observe that for a regular Bratteli diagram with an ordering $\om$, the Vershik map $\varphi_\om$, if it exists, is
defined in a unique way. Also, a necessary condition for $\om\in
\mathcal P_{B}$ is that $|X_{\max}(\om)|=|X_{\min}(\om)|$. Given $(B,\om)$ with $\om \in \mathcal P_B$, the uniquely defined  system $(X_B, \varphi_\om)$ is called  a {\em Bratteli-Vershik} or {\em adic} system.

Thus, we can summarize the above definitions and results in the following  statement.

\begin{theorem}  Let $B= (V,E, \om)$ be an ordered Bratteli diagram with perfect order
$\om  \in \mathcal P_B$. Then there exists an aperiodic homeomorphism (Vershik map) $\varphi_{\om}$ acting on the path space $X_B$ according to Definition \ref{VershikMap}. The homeomorphism $\varphi_\om$ is minimal if and only if $B$ is simple.
\end{theorem}

The pair $(X_B, \varphi_\om)$ is called the {\em Bratteli-Vershik dynamical system}.

The simplest example of a Bratteli diagram is an {\em odometer}. Any odometer can be realized as a Bratteli diagram $B$ with $|V_n| =1$ for all $n$. Then any order on $B$ is proper and defines the Vershik map.

It is worth noticing that a general Bratteli diagram may have a rather complicated structure. In particular, the tail equivalence relation may have uncountably many minimal components or, in other words, uncountably many simple subdiagrams that do not have connecting edges.

\section{Bratteli-Vershik representations of Cantor dynamical systems
and orbit equivalence}\label{BV}

\subsection{Bratteli-Vershik model of a Cantor dynamical system}
Let $(X,T)$ be a Cantor aperiodic dynamical system. Is it possible to represent $(X,T)$ as a Bratteli-Vershik system? In other words, we want to associate a refining sequence of clopen partitions to the homeomorphism $T$ whose elements generate the clopen topology on $X$. If $T$ is minimal, the answer is well known. Take any  clopen subset $A$ of $X$ and consider the forward $T$-orbit that starts at a point $x \in A$. By minimality of $T$, the orbit returns to $A$ in a finitely many steps. Since the function of the first return $r_A(x)$ has finite values and $A$ is compact, we get a finite partition of $X$ into clopen $T$-towers $\xi_i = (A_i, TA_i, ... , T^{i-1}A)$ with base $A = \bigcup_i A_i$ where $A_i = \{x\in A : r_A(x) = i\}$.  For aperiodic $T$, this partition can be obtained only when  $A$ has some additional properties.

We say that $\xi$ is a {\em Kakutani-Rokhlin partition} of a Cantor set $X$ if $\xi$ is a finite union of disjoint $T$-towers $\xi_i = (C_i, ..., T^{h_i -1}C_i)$.  Then $h_i$ is called the height of $\xi$ and $C_i$ is the base of $\xi_i$.

One says that a partition $\xi_1$ refines a partition $\xi$ if every element (atom) of $\xi$ is a union of elements of $\xi_1$.

A closed subset $Y$ of $X$ is called {\em basic} if (1) $Y \cap T^iY = \emptyset, i\neq 0$, and (2) every clopen neighborhood $A$ of $Y$ is a complete $T$-section, i.e., $A$ meets every $T$-orbit at least once. This means that every point from $A$ is recurrent. It is clear that if $T$ is minimal then every point of $X$ is a basic set.

The following result was well known for a minimal homeomorphism $T$ after the paper by Putnam \cite{putnam:1989}. The case of aperiodic Cantor system is much subtler and was considered in \cite{BDM05}  and \cite{medynets:2006}. Based on the results of \cite{BDM05} Medynets proved in \cite{medynets:2006} that {\em every Cantor aperiodic system $(X, T)$ has a basic set}. Then one can prove

\begin{theorem} [\cite{medynets:2006}] \label{K-R part} Let $(X,T)$ be an aperiodic Cantor system, and let $Y \subset X$ be a basic set for $T$. Then there exists a sequence of Kakutani-Rokhlin clopen partitions $\xi(n)$ such that for all $n \in \N$:

(i) $\xi(n+1)$ refines $\xi(n)$ and atoms of these partitions generate the clopen topology on $X$;

(ii) $B(\xi(n)) \subset B(\xi(n+1)) $ and $Y = \bigcap_n B(\xi(n))$;

(iii) $h_n \to \infty$ (as $n\to \infty$) where $h_n$ is the minimal height of the towers that form $\xi(n)$.
\end{theorem}

The ideas developed in the papers by Vershik \cite{vershik:1981}, \cite{vershik:1982}, where sequences of refining measurable partitions of a measure space were used to construct a realization of an ergodic automorphisms of a measure space, turned out to be very fruitful for finding a model of any minimal homeomorphism $T$ of a Cantor set $X$. In \cite{herman_putnam_skau:1992}, Herman, Putnam, and Skau found an explicit construction that allows one to define an ordered simple Bratteli diagram $B = (V, E, \om)$ such that $T$ is conjugate to the corresponding Vershik map $\varphi_\omega$. Since this construction is described in many papers (not only in \cite{herman_putnam_skau:1992}), we will not  give the details here referring to the original paper \cite{herman_putnam_skau:1992} (and \cite{durand:2010}) for detailed explanation. We discuss here the main idea of constructing such a diagram $B$ in a non-rigorous way.

Let $T$ be an aperiodic (minimal) homeomorphism of a Cantor set $X$. First,  find a sequence of Kakutani-Rokhlin partitions $(\xi(n))$ satisfying Theorem \ref{K-R part}: $\xi(0) = X$,
 \begin{equation*}\xi(n)=\{T^iA_v(n)\; :\; v=1,\ldots, m_v(n); \ \  i=0,\ldots,h_v(n)-1\},\;n\geq 1,
\end{equation*}
that generates the clopen topology on $X$. Simultaneously, we will define an ordered  Bratteli diagram $B = (V,E, \om)$ as follows.

(i) Let
$$
\xi_v(n)=\{A_v(n),\ldots,T^{h_v(n)-1}A_v(n)\} \mbox{ for } v=1,\ldots, m_v(n).
$$
The vertex set $V = \coprod_{n \geq 0} V_n$ where $V_0$ is a singleton, and  $V_n  = \{1, ... , m_v(n)\}$.

(ii) Define the set  of edges $E_n$  between the consecutive
levels $V_{n}$ and $V_{n+1}$ by the incidence matrix
$F_n=\{m_{vw}(n) : v\in V_{n+1},\;w\in V_{n}\}$, where
$$
f_{vw}^{(n)}= |\{0\leq i<h_v(n) : T^iA_v(n+1)\subset A_w(n)\}|.
$$
In other words, we fix a vertex $v\in V_{n+1}$ and define $V(v,n)$
as the set of all vertices from $V_{n}$ such that $\xi_v(n +1)$
intersects $\xi_w(n)$, and a vertex $w$ appears in $V(v,n)$ as
many times as $\xi_v(n+1)$ intersects $\xi_w(n)$. Then we connect
$v$ to each vertex $w\in V(v,n)$ taking into account the
multiplicity of appearance of $w$ in $V(v,n)$.

(iii) To define the ordering $\om$ on $E$ we take the clopen set
$A_v(n+1)$. Then tracing the orbit of $A_v(n+1)$ within the $T$-tower
$\xi_v(n+1)$, we see that $T^iA_v(n+1)$ consecutively meets the sets
$A_{w_1}(n)$, ..., $A_{w_{k_v}}(n)$ (some of them can
occur several times). This defines the set of edges $r^{-1}(v)$.
Enumerate the edges from $r^{-1}(v)$ as follows:
$e(w_1,v) < e(w_2,v) <\ldots<e(w_{k_v},v)$.

Since the partitions $\{\xi_n\}$ generate the topology of
$X$, for each point $x\in X$ there is a unique sequence
$i(x)=\{(v_n,i_n)\}_{v\in V_n;\,0\leq i_n<h_{v_n}(n)}$ such that
$$
\{x\}=\bigcap_{n=1}^\infty T^{i_n}A_{v_n}(n).
$$
Define the map $\theta : X\rightarrow X_B$ by
$$\theta(x)=\bigcap_{n\geq 1}U(y_1,\ldots,y_n)
$$ where
$(y_1,\ldots,y_n)$ is the $i_n$-th finite path in $E(v_0,v_n)$
with respect to the lexicographical ordering on $E(v_0,v_n)$.  It can be easily checked that $\theta$ is a homeomorphism.

Define $\varphi_\om=\theta \circ T\circ\theta^{-1}$. Then $\varphi_\om$ is the
Vershik map defined by order $\om$. Thus, we
obtain that $(X,T)$ is conjugate  to the Bratteli-Vershik system
$(X_B,\varphi_\om)$.

This gives the theorem proved in \cite{medynets:2006} for an aperiodic homeomorphism (the case of a minimal homeomorphism was considered in  \cite{herman_putnam_skau:1992}).

\begin{theorem} \label{existenceBD}
Let $(X, T)$ be a Cantor aperiodic system with a basic set $Y$. There exists an ordered Bratteli diagram  $B = (V,E,\om)$ such that $(X, T)$ is conjugate to a Bratteli-Vershik dynamical system $(X_B, \varphi_\om)$. The homeomorphism implementing the conjugacy between $T$ and $\varphi_\om$ maps the basic set $Y$ onto the set $X_{\min}(\om)$ of all minimal paths of $X_B$. The equivalence class of $B$ does not depend on a choice of $\{\xi(n)\}$ with the property $\bigcap_n B(\xi(n)) =Y$.

\end{theorem}

Thus, for every aperiodic homeomorphism $T$ of a Cantor set $X$, there exists an ordered Bratteli diagram $(B, \om)$ such that $T$ is conjugate to the Vershik map $\varphi_\om$. Is the converse true? In case of a simple Bratteli diagram, the answer is  obviously affirmative: there exists a proper order $\om$ on any simple Bratteli diagram $B$ so that $(X_B, \varphi_\om)$ is a minimal Cantor system.

The main difficulty for the study of non-simple Bratteli diagrams is illustrated by the following amazing observation made by Medynets in \cite{medynets:2006}.
He proved that there exists an (even stationary) non-simple Bratteli diagram such the {\em Vershik map $\varphi_\om$ is not continuous for any order $\om$.} Here is an example of such a diagram:

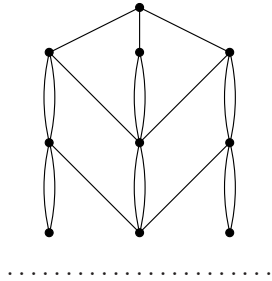
\begin{center}
\unitlength=0.6cm
\begin{graph}(5,7)
 \graphnodesize{0.2}
 \roundnode{V0}(3,6)
 \roundnode{V11}(1,5)
 \roundnode{V12}(3,5)
 \roundnode{V13}(5,5)

 \roundnode{V21}(1,3)
 \roundnode{V22}(3,3)
 \roundnode{V23}(5,3)

 \roundnode{V31}(1,1)
 \roundnode{V32}(3,1)
 \roundnode{V33}(5,1)

 %
 \graphlinewidth{0.025}
 \edge{V0}{V11}
  \edge{V0}{V12}
  \edge{V0}{V13}

 \bow{V21}{V11}{0.06}
  \bow{V21}{V11}{-0.06}
 \edge{V22}{V11}
 \bow{V22}{V12}{0.06}
 \bow{V22}{V12}{-0.06}
 \edge{V22}{V13}
 \bow{V23}{V13}{0.06}
 \bow{V23}{V13}{-0.06}

 \bow{V31}{V21}{0.06}
  \bow{V31}{V21}{-0.06}
 \edge{V32}{V21}
 \bow{V32}{V22}{0.06}
 \bow{V32}{V22}{-0.06}
 \edge{V32}{V23}
 \bow{V33}{V23}{0.06}
 \bow{V33}{V23}{-0.06}
 \freetext(3,0.1){.\,.\,.\,.\,.\,.\,.\,.\,.\,.\,.\,.\,.\,.\,.\,.\,.\,.\,.\,.\,.\,.\,.}

\end{graph}

\end{center}

\subsection{Orbit equivalence and full groups}

The circle of problems related to the classification of transformations with respect to orbit equivalence was originated in the famous papers by Dye \cite{Dye1}, \cite{Dye2}. It was proved that any two ergodic finite measure preserving automorphisms of a standard measure space are orbit equivalent. Also he defined a new invariant of orbit equivalence, the so called {\em full group} of automorphisms. We give its definition for homeomorphisms of a Cantor set, but it is obvious that this definition can be reformulated for dynamical systems acting on other underlying spaces as a Borel space or standard measure space.

\begin{definition}\label{Full GR def}
Let $\Gamma$ be  a group of homeomorphisms of a Cantor set $X$. Then the full group $[G]$ (generated by $\Gamma$) is formed by homeomorphisms of $X$ that keep the partition of $X$ into $\Gamma$-orbits fixed.  That is
$$
[\Gamma] := \{R \in Homeo (X) : Rx \in Orb_{\Gamma}(x), \forall x \in X\}.
$$
\end{definition}

The most interesting is the case when $\Gamma$ is an action of an amenable countable group $G$. If $\Gamma$ is generated by a single homeomorphism $T$ we write $[T]$ for the corresponding full group. It is well known that in the case of measurable dynamics any ergodic finite measure preserving $G$-action is orbit equivalent to an action of $\Z$ \cite{OW80}, \cite{CFW81}. One of the central problems of Cantor dynamics consists of an extension of this result to the group of homeomorphisms of a Cantor set.

In 1995, Giordano, Putnam, and Skau \cite{giordano_putnam_skau:1995} found complete invariants of orbit equivalence for minimal homeomorphisms. We need to give more definitions to formulate their results. Giving a Cantor minimal system $(X, T)$, denote by $C(X, \Z)$ the countable group of integer-valued  continuous functions, and by $B(X,\mathbb Z)$ the subgroup of $C(X,\Z)$ generated by $T$-coboundaries ($f \in B(X,\mathbb Z)$ if $f = g - g\circ T$ for some $g\in C(X,\mathbb Z)$).
Then
$$
K^0(X, T) = C(X, \mathbb Z)/B(X,\mathbb Z)
$$
is considered as an ordered abelian group whose cone of positive elements $K^0(X, T)^+$ is determined by $\{ \widehat f \in K^0(X, T) : f \geq 0, f \in C(X, \mathbb Z)\}$ with distinguished order unit $\widehat 1$. The group $C(X, \Z)$ contains also the subgroup of all infinitesimal functions with respect to $T$: by definition, $f \in Inf(X, T)$ if
$$
\int_X f d\mu = 0 \ \mbox{ for\ all} \ \mu\in M(X,T).
$$
Since $Inf(X, T) \supset B(X,T)$, one can define the ordered abelian group
$$
\wt K^0(X, T) = K^0(X, \mathbb Z)/Inf(X,\mathbb Z).
$$

\begin{theorem}[\cite{giordano_putnam_skau:1995}, \cite{GPS99}] \label{OE GPS}
Let $(X,T)$ and $(Y,S)$ be minimal Cantor systems. The following are equivalent:

(1) $(X,T)$ and $(Y,S)$ are orbit equivalent;

(2) the groups $\wt K^0(X, T)$ and $\wt K^0(Y, S)$ are order isomorphic by a map preserving the distinguished order units;

(3) there exists a homeomorphism $F : X \to Y$ carrying the $T$-invariant measures onto the $S$-invariant measures;

(4) the full groups $[T]$ and $[S]$ are isomorphic as abstract groups.
\end{theorem}

There are two impressive corollaries of this theorem related to uniquely ergodic homeomorphisms. Let $S(\mu) := \{\mu (A) : A \ \mbox{clopen}\}$ where $\mu$ is a probability measure on a Cantor set $X$.

\begin{corollary} [\cite{giordano_putnam_skau:1995}] \label{OE corollary}
(1) The uniquely ergodic Cantor minimal systems $(X,T)$ and $(Y,S)$ are orbit equivalent if and only if $S(\mu) = S(\nu)$ where $ \mu \circ T = \mu$ and $\nu\circ T = \nu$.

(2) Any uniquely ergodic minimal Cantor system is orbit equivalent either an odometer (when $S(\mu) \subset \mathbb Q$) or to a Denjoy  homeomorphism.
\end{corollary}

An action $\Gamma$ of a group $G$ is called {\em hyperfinite or affable (approximately finite)} if it is orbit equivalent to a $\Z$-action.  It is an intriguing  question which countable groups admit hyperfinite actions by homeomorphisms of a Cantor set. So far, the best known result is

\begin{theorem} [\cite{GMPS10}] \label{Z^n} Any minimal action of $\Z^n$ on a Cantor set is hyperfinite.
\end{theorem}

Hyperfinite equivalence relations have been studied in a series of papers that contain very deep results on their structure like the Absorbtion Theorem (see \cite{GMPS08}). We refer to a recent paper by Putnam \cite{P10} devoted to this subject.

The paper \cite{HKY12} is devoted to the interplay between the notions of topological and measurable orbit equivalence.  Let $(X,T)$ and $(Y,S)$ be two Cantor minimal systems. A homeomorphism $\rho\:\ X\rightarrow Y$ is called universally measure preserving if $\rho$ gives rise to an affine bijection between $T$-invariant probability measures and $S$-invariant probability measures.   The main result of the paper is a dynamical proof of the fact that two Cantor minimal systems $(X,T)$ and $(Y,S)$ are orbit equivalent if and only if there exists a universally measure preserving homeomorphism between $X$ and $Y$. The proof is based on a scrupulous Kakutani-Rokhlin tower analysis (in its Bratteli-Vershik interpretation). Another important ingredient of the proof is a topological version of the famous copying lemma  which is developed in \cite{HKY12}. 
\medskip

It is unclear to what extend this result can be generalized. This is an open problem whether any minimal action of a countable amenable group is hyperfinite. Moreover very few facts are known about orbit equivalence of aperiodic homeomorphisms of a Cantor set. We will see below that the Theorem \ref{OE GPS} is not true for aperiodic homeomorphisms. It is obvious that aperiodic actions have more invariants of orbit equivalence than minimal actions have. For example, the cardinality of the set of minimal components must be the same for orbit equivalent aperiodic homeomorphisms.

To finish this subsection, we mention  the following result proved by Medynets in \cite{Med}. It is proved  that (under some very mild conditions) the Cantor dynamical systems $(X, G_1)$ and $(X, G_2)$ are orbit equivalent if and only if the full groups $[G_1]$ and $[G_2]$ are isomorphic as abstract groups. This theorem, in particular, works for aperiodic actions of countable groups.

\section{Bratteli diagrams that admit a continuous Vershik map}\label{VM}

In this section we are interested in the following problem. How can we recognize whether a given Bratteli diagram admit an order generating a continuous Vershik map? In other words, we are interested in the problem of finding those Bratteli diagrams that can serve as Bratteli-Vershik models of aperiodic homeomorphisms.

We can  describe the set of all orderings $\mathcal O_{B}$ on a Bratteli diagram $B$ in the following way. Let   $P_v$ denote the set of all orders on $r^{-1}(v)$; an element in $P_v$ is denoted by $\om_v$. Then $\mathcal O_{B}$ can be
represented as
\begin{equation}\label{orderings_set}
\mathcal O_{B} = \prod_{v\in V}P_v .
\end{equation}
 Giving each set $P_v$ the discrete topology, we see from (\ref{orderings_set}) that  $\mathcal O_{B}$ is a Cantor set with respect to the product topology. In other words, two orderings $\omega= (\om_v)$ and $\omega' = (\om'_v)$ from $\mathcal O_{B}$ are close
if and only if they agree on a sufficiently long initial segment: $\om_v = \om'_v,
v \in \bigcup_{i=0}^k V_i$.

It is worth  noticing that the order space $\mathcal O_B$ is sensitive with
respect to a telescoping. Indeed, let $B$ be a Bratteli diagram and $B'$
denote the diagram obtained
by telescoping of $B$ with respect to a subsequence $(n_k)$ of levels.
We see that any ordering $\omega$ on $B$ can be extended to the
(lexicographic) ordering $\om'$ on $B'$. Hence the map $L : \om \to \om'=L(\om)$
defines a closed proper subset $L(\mathcal O_B)$ of $\mathcal O_{B'}$.

The set of all orderings $\mathcal O_{B} $ on a Bratteli diagram $B $
can be considered also as a \textit{measure space} whose Borel structure is
generated by cylinder sets.  On the set
$\mathcal O_B $
we take the product measure $\mu = \prod_{v\in V} \mu_v$ where  $\mu_v$
is a measure on the set $P_v$.   The case where each $\mu_v$
is the uniformly distributed measure on $P_v$ is of particular interest:
$\mu_v(\{i\}) =  (|r^{-1}(v)|!)^{-1}$ for every $i \in P_v$ and $v\in V^* \backslash V_0$. Unless $|V_n|=1$ for almost all $n$,  if $B'$ is a telescoping of $B$, then in $\mathcal O_{B'}$, $L(\mathcal O_{B})$ is a set of zero measure.

We recall that $\mathcal P_B$ is the  subset of $\mathcal O_B$ consisting of \textit{perfect} orderings that produce Bratteli-Vershik topological dynamical systems (Vershik maps).

One can show that the map $\Phi : \omega \to \varphi_\om $ is continuous.
Also the sets $\mathcal P_B$ and ${\mathcal O_B} \setminus {\mathcal P}_B$ are both dense in $\mathcal O_B$ (see \cite{BKY14}).

In \cite{BKY14}, a class of Bratteli diagrams that do not have perfect orders was described.

\begin{definition}\label{aperiodic_definition}
We define the family $\mathcal A$ of Bratteli diagrams, all of whose incidence matrices are of the form
\[F_{n}:=
\left(
\begin{array}{ccccc}
A_{n}^{(1)} & 0 & \ldots  &  0 & 0 \\
0 & A_{n}^{(2)} & \ldots  & 0 & 0 \\
\vdots & \vdots  & \ddots & \vdots & \vdots     \\
0 & 0 & \ldots  & A_{n}^{(k)}  & 0 \\
B_{n}^{(1)}& B_{n}^{(2)}  & \ldots  & B_{n}^{(k)} & C_{n}  \\
 \end{array}
\right), \ \ n\geq 1,
\]
where
\begin{enumerate}
\item
for $1\leq i \leq k$ there is some $d_{i}$ such that
$A_{n}^{(i)}$ is a $d_{i}\times d_{i}$ matrix for each $n\geq 1$,
\item all  matrices $A_{n}^{(i)}$, $B_{n}^{(i)}$ and $C_{n}$ are strictly
 positive,
\item $C_{n}$ is a $d\times d$ matrix,
\item there exists $j\in \{\sum_{i=1}^{k}d_{i} +1, \ldots \sum_{i=1}^{k}d_{i} +d\}$ such that for each $n\geq 1$, the $j$-th row of $F_{n}$ is strictly positive.
\end{enumerate}
If a Bratteli diagram's incidence matrices are of the form above, we shall say that it has $k$ {\em minimal components}.
\end{definition}

It follows from this definition that every Bratteli diagram has exactly $k$ minimal components.

The next proposition describes how for some aperiodic diagrams $B$ that
belong to the special class $\mathcal A$ (see Definition \ref{aperiodic_definition}), there are structural obstacles to the existence of perfect orders on $B$.

\begin{proposition}[\cite{BKY14}]\label{corollary_1}
Let $B\in \mathcal A$ have $k$ minimal components, and such that for each $n\geq 1$,  $C_{n}$ is an $s\times s$ matrix where $1\leq s\leq k-1$. If $k>2$, then there is no perfect ordering on $B$. If $k=2$, there are perfect orderings on $B$ only if $C_{n}=(1)$ for all but finitely many $n$.
\end{proposition}

The technique used in \cite{BKY14} and \cite{BY} is based on study new notions related to any ordered Bratteli diagram. They are  skeletons and associated graphs.
These notions are especially useful in the case of finite rank diagrams.

Suppose that  $B$ has strict rank $d$, i.e., $|V_n| = d$ for all $n\geq 1$.
  If a maximal (minimal) path $M$ ($m$) goes through the same vertex $v_{M}$
($v_{m}$) at each level of $B$, we will call this path {\em
vertical}. The following proposition characterizes when $\om$ is
a perfect order on a finite rank Bratteli diagram.

\begin{proposition} [\cite{BKY14}] \label{existence_vershik_map}
Let $(B, \om)$ be an ordered Bratteli diagram.

\begin{enumerate}
\item Suppose that $B$ has strict rank $d$ and that the
  $\om$-maximal  and $\om$-minimal paths  $M_1,...,M_k$ and
      $m_1,...,m_{k'}$
are vertical passing through the vertices $v_{M_1},\ldots , v_{M_k}$ and $v_{m_1},\ldots , v_{m_{k'}}$ respectively.  Then $\om$ is perfect if and only if

\begin{enumerate}\item $ k = k'$,

\item there is a permutation $\sigma$ of $\{1,\ldots k\}$ such that for
each $i\in \{1,...,k\}$, $v_{M_i}v_{m_{j}} \, \in \mathcal L(B,\om)$ if and only if $j=\sigma(i)$.

\end{enumerate}
\item
  Let  $B'$ be  a telescoping of
 $B$. Then $\omega \in \mathcal P_{B}$ if and only if $\om' = L(\om)\in \mathcal P_{B'}$.

\end{enumerate}
\end{proposition}

Let  $\om$ be an order on a  Bratteli diagram $B$.    If  $v\in V\backslash V_0  $, we denote the  minimal edge with range $v$ by  $\ol e_v$ , and  we denote the maximal edge with range $v$ by $\wt e_v$.

\begin{definition} \label{well_telescoped_definition}
Let $(B,\om)$ be an ordered  rank $d$ diagram. We say that $(B,\om)$ is {\em well telescoped} if
\begin{enumerate}

\item $B$ has strict rank $d$,

\item all $\om$-extremal paths are vertical, with $\wt V$,  $\ol V$ denoting the sets of vertices through which maximal and minimal paths run respectively,

\item $s(\wt e_v)\in \wt V$ and $s(\ol e_v) \in \ol V$ for each $v\in V\backslash (V_0 \cup V_1 )$, and this is independent of $n$.
\end{enumerate}
If $(B,\om)$ is perfectly ordered, for it to be  considered well telescoped, it will also have to satisfy
\begin{enumerate}
\setcounter{enumi}{3}
\item if $\wt v\ol v $ appears as a subword of some $w(v,m,n)$ with $m\geq 1$, then , then $\sigma(\wt v) = \ol v$ defines a one-to-one correspondence between the sets $\wt V$ and $\ol V$.
\end{enumerate}

\end{definition}
Given an ordered finite rank $(B,\om)$, it can always be telescoped so that it is well telescoped.
For details of how this can be done, see Lemma 3.11 in \cite{BKY14}. Thus, when we talk about a (finite rank) ordered diagram, we assume without loss of generality  that it is well telescoped.
For well telescoped ordered diagrams $(B,\om)$, we have $s(\wt e_v) \in \wt V_n$ and $s(\ol e_v) \in \ol V_n$ for any $v\in V_{n+1},\ n \geq 1$.

\begin{definition}
Given a well telescoped $(B,\om)$, we call the set
$$
\mathcal F_{\omega} = (\wt V, \ol V, \{\wt e_{v}, \ol e_{v} : v\in V_{n}, \,\, n\geq 2\})
$$
the {\em skeleton associated to $\om$}. If  $\om$ is a perfect order on $B$,   it follows that $|\ol V |= |\wt V|$, and
if $\sigma:\wt V\rightarrow \ol V$ is the
permutation given by Proposition \ref{existence_vershik_map}, we call
$\sigma$ the {\em accompanying permutation}.
\end{definition}

The notion of a skeleton of an ordered diagram can be extended to an unordered diagram. Namely, given a strict rank $d$ diagram $B$, we select, two subsets $\wt V$ and $\ol V$ of $V$, of the same cardinality, and, for each $v\in V\backslash V_0 \cup V_1$, we select two edges $\wt e_v$ and $\ol e_v$, both with range $v$, and such that  $s(\wt e_v) \in \wt V$,
$s(\ol e_v) \in \ol V$. In this way we can extend the definition of a skeleton
$\mathcal F = (\wt V, \ol V, \{\wt e_{v}, \ol e_{v}: v\in V_{n}, \,\, n\geq 2\})$ to an unordered strict rank Bratteli diagram.
 A  more detailed discussion  can be found in \cite{BKY14}).  Arbitrarily choosing a
bijection $\sigma:\wt V \rightarrow \ol V$, we can consider the set of orders on $B$ which have $\mathcal F$ as skeleton and $\sigma$ as accompanying permutation.

\begin{definition} Given a skeleton $\mathcal F$ on a strict finite rank diagram $B$,
for any vertices $\wt v \in
\wt V$ and $\ol v \in \ol V$, we set
\begin{equation}\label{W} W_{\wt v} = \{w \in V : s(\wt e_w) =
\wt v\} \mbox{ and }  W'_{\ol v} = \{w \in V : s(\ol e_w) =
\ol v\}.
\end{equation}
 Then $W = \{W_{\wt v} : \wt v \in \wt V\}$ and $W' = \{W_{\wt v}' : \ol v \in \ol V\}$ are both partitions of $V$.  We call $W$ and $W'$ the {\em partitions generated by $\mathcal F$.}
 \end{definition}

 Let $[\ol v, \wt v]:= W'_{\ol v} \cap W_{\wt v}$, and
define the partition \[W \cap W' := \{[\ol v,\wt v] : \ol v\in \ol V, \wt v \in \wt V, [\ol v,\wt v]\neq \emptyset\} .\]

\begin{definition}\label{associated_directed_graph}
Let $\mathcal F$  be a skeleton on
the strict finite rank  $B$ with
accompanying permutation $\sigma$.
Let  $\mathcal H=(T,P)$ be the directed graph where
the set $T$ of vertices of  $\mathcal H$  consists of partition elements $[\ol v, \wt v]$ of $W'\cap W$, and where there is an edge in  $P$ from $[\ol v, \wt v]$ to $[\ol v', \wt v']$
if and only if $\ol v' \, =\sigma(\wt v)$.
We call $\mathcal H$ the {\em  directed graph associated to $(B, \mathcal F, \sigma)$}.
\end{definition}

In order to see how these notion work, we can formulate the following illustrating results. Recall that a directed graph is {\em strongly connected} if for any two
vertices $v$, $v'$, there is a  path from $v$ to $v'$, and also  a path from $v'$
to $v$. If at least one of these paths exist, then $G$ is {\em weakly
connected}.  We notice that, given $(B, \mathcal F, \sigma)$,  an associated graph $\mathcal H =(T,P)$ is not connected, in general.

\begin{proposition}[\cite{BKY14}]\label{connectedness and goodness}
 Let $(B,\om)$ be a finite rank,  perfectly ordered and well telescoped  Bratteli diagram, and suppose $\om$ has skeleton $\mathcal F_\om$ and permutation $\sigma$.
\begin{enumerate}
\item If $B$ is simple, then the associated graph $\mathcal H$ is
strongly connected.
\item If $B \in \mathcal A$, then the associated graph $\mathcal H$ is
weakly connected.
\end{enumerate}
\end{proposition}

If we wanted to define similar notions for a not finite rank Bratteli diagram $B$, then we would come across the obvious difficulties because $|V_n|$ is not bounded. Nevertheless, one can overcome these problems and define a skeleton of $B$ and a sequence of associated graphs. The details are rather cumbersome, so that we refer to \cite{BY} where these definitions and more facts were given.

Suppose a Bratteli diagram $B$  is given. The structure of the diagram is completely defined by the sequence of incidence matrices $(F_n)$. The following question seems to be challenging: {\em is it possible to determine using the entries of $F_n$  whether the Bratteli diagram $B$ admits a perfect order?} We give a criterion that answers the above question based on \cite{BKY14} (for a finite rank diagram) and \cite{BY} (for an arbitrary Bratteli diagram). Here we discuss the case of a finite rank diagram since it is more transparent, and the general result from \cite{BY}  extends the method used in \cite{BKY14}.

Let $(B,\om)$ be a perfectly ordered simple and well telescoped Bratteli diagram of finite rank diagram (that is $|V_n|  = d$). Let  $\mathcal F = \mathcal F_\om$ be the  skeleton  generated by $\om$,  $\sigma:\wt V \rightarrow \ol V$ the permutation, and $\varphi = \varphi_\om$ be the corresponding Vershik map such that $\varphi_{\om} (M_v) = m_{\sigma(v)}$ for $ v \in \wt V$.
Also we define the two partitions $W = \{W_{\wt v} : v \in \wt V\}$ and $W' = \{W'_{\ol v} : v \in \ol V\}$ of $V$ generated by $\mathcal F$.

 Let $E(V_n, u)$ be the set of all finite paths
between vertices of level $n$ and a vertex $u\in V_m$ where $m >
n$. The symbols $\wt e(V_n, u)$ and $\ol e(V_n, u)$ are used to denote
the maximal and minimal finite paths in $E(V_n, u)$, respectively.
Fix maximal and minimal vertices $\wt v$ and $\ol v$ in $\wt V_{n-1}$ and $\ol V_{n-1}$ respectively. Denote $E(W_{\wt
v}, u) = \{e \in E(V_n, u): s(e) \in W_{\wt v}, r(e) = u\}$ and
$\wt E(W_{\wt v}, u) = E(W_{\wt v}, u) \setminus \{\wt e(V_n,
u)\}$. Similarly, $\ol E(W'_{\ol v}, u) = E(W'_{\ol v}, u)
\setminus \{\ol e(V_n, u)\}$. Clearly, the sets $\{E(W_{\wt v}, u)
: \wt v \in \wt V\}$ form a partition of $E(V_n, u)$.
Let $e$ be a non-maximal finite path, with  $r(e)=v$ and $s(e)\in V_m$,  which determines the cylinder set $U(e)$.  It is clear that for any finite path $e \in \wt E(W_{\wt v}, u)$ we have $\varphi_\om(e) \in \ol E(W'_{\sigma(\wt v)}, u)$.
Therefore,
$$
|\wt E(W_{\wt v}, u)| = |\ol E(W'_{\sigma(\wt v)}, u)|.
$$

Define two sequences of matrices $\wt F_n = (\wt
f_{w,v}^{(n)})$ and $\ol F_n = (\ol f_{w,v}^{(n)})$ by the following
rule (here $w\in V_{n+1}$, $ v\in V_n$ and $ \ n\geq 1$):
\begin{equation}\label{wt f}
\wt f_{w,v}^{(n)} = \left\{
                      \begin{array}{ll}
                        f_{w,v}^{(n)} - 1, & \mbox{ if }\wt e_w \in E(v, w); \\
                        f_{w,v}^{(n)}, & \mbox{otherwise},
                      \end{array}
                    \right.
\end{equation}
\begin{equation}\label{ol f}
\ol f_{w,v}^{(n)} = \left\{
                      \begin{array}{ll}
                        f_{w,v}^{(n)} - 1, & \mbox{ if }\ol e_w \in E(v, w); \\
                        f_{w,v}^{(n)}, & \mbox{otherwise}.
                      \end{array}
                    \right.
\end{equation}

Then for any $u\in V_{n+1}$ and $ \wt v\in \wt V_{n-1}$, we obtain that the entries of incidence matrices have the property:
\begin{equation}\label{sum condition1}
\sum_{w\in W_{\wt v}} \wt f_{u,w}^{(n)} \ \  = \sum_{w'\in W'_{\sigma(\wt v)}} \ol f_{u,w'}^{(n)},\ \ n \geq 2.
\end{equation}
We call  relations (\ref{sum condition1}) the {\em balance relations}.

\begin{theorem} [\cite{BKY14}] \label{existence of good order}
Let $B$ be a simple strict rank $d$ Bratteli diagram,  let $\mathcal F = \{M_{\wt v}, m_{\ol v}, \wt e_w, \ol e_w: w \in V \backslash V_0,\,\, \wt v \in \wt V\mbox{  and } \ol v \in \ol V  \}$ be a skeleton on $B$, and let $\sigma : \wt V \to \ol V$ be a bijection.  Suppose that eventually all associated graphs $\mathcal H_n$ are positively strongly connected, and suppose that the entries of incidence matrices $(F_n)$ eventually satisfy the balance relations (\ref{sum condition1}).
Then there is a perfect ordering $\om$ on $B$ such that $\mathcal F = \mathcal F_\om$ and the Vershik map $\varphi_\om$ satisfies the relation $\varphi_\om(M_{\wt v}) = m_{\sigma(\wt v)}$.
\end{theorem}

The next results is related to finite rank Bratteli diagrams. Recall that on the set of all orders $\mathcal O_B$ we can consider the product measure $P$ such that the probability to pick up an order for $r^{-1} (v) $ is uniform. The following theorem is somewhat surprising because it states that the diagram $B$ ``knows a priori'' how many maximal and minimal paths it must have.

\begin{theorem} [\cite{BKY14}]\label{goodbaddichotomy}
 Let $B$ be a finite rank $d$ aperiodic Bratteli diagram. Then there
 exists  $j \in \{1,..., d\}$ such that $P$-almost all orders
 have $j$ maximal and $j$ minimal elements.
\end{theorem}

In the recent work \cite{JQY14}, random orders have been studied on simple 0-1 Bratteli diagrams, i.e., all entries of incidence matrices are either zeros or ones. 
 It was proved that a random order has uncountably many infinite paths if and only if the growth rate of the level-n vertex sets is superlinear.

\begin{theorem}[\cite{JQY14}] \label{JQY} Let $B$ be a Bratteli diagram with incidence matrices $F_n$ whose entries are 1 and  $|V_n| \geq  1$. The following dichotomy holds:

(1) If $\sum_n 1/|V_n| = \infty$, then there exists  $P$-almost sure a unique maximal path.

(2) If $\sum_n 1/|V_n| < \infty$, then there exists  $P$-almost sure
uncountably many maximal paths.

\end{theorem}

Thus, a random order on a slowly growing Bratteli diagram satisfying \ref{JQY}
admits a homeomorphism, while a random order on a quickly growing Bratteli diagram does not.

\section{Stationary Bratteli diagrams}\label{stationary}

In this section, stationary Bratteli diagrams are considered. We explicitly describe the class of homeomorphisms represented by stationary simple and non-simple diagrams. We also consider their invariant measures and orbit equivalence classes.

\subsection{Substitution dynamical systems}

Papers~\cite{forrest, durand_host_skau} explore minimal Cantor systems and their relation to stationary Bratteli diagrams and the corresponding dimension groups (for example, see~\cite{Effros} for details on dimension groups). It turns out that the class of minimal homeomorphisms which can be represented by stationary Bratteli diagrams is constituted by minimal substitution dynamical systems and odometers. In \cite{BKM09}, the systematical study of non-primitive substitutions was initiated. In particular, it is shown that one can prove the analogue of the above mentioned result for aperiodic homeomorphisms, i.e., the assumption of minimality is not essential for finding the corresponding Bratteli-Vershik model.

Let $A$ denote a finite alphabet and $A^+$ the set of all non-empty words over $A$. For a word $w=w_0\ldots w_{n-1}$
with $w_i\in A$ let $|w|=n$ stand for its length. By a {\it substitution}, we mean any map $\sigma : A\rightarrow A^+$. A map $\sigma : A\rightarrow A^+$ can be extended to the map $\sigma :A^+\rightarrow A^+$ by concatenation. We define the {\it language} $L(\sigma)$ of  a substitution $\sigma$ as the set of all words which appear as factors of $\sigma^n(a)$, $a\in A$, $n\geq 1$. We set also $\sigma^0(a)=a$ for all $a\in A$. Let $T : A^\mathbb {Z} \to A^\mathbb {Z}$ be the shift, that is $T(x_i) = (y_i)$ where $y_i = x_{i+1}$, $i\in \Z$. By a {\it substitution dynamical system associated to a substitution $\sigma$}, we mean a pair $(X_\sigma,T_\sigma)$, where
$$
X_\sigma=\{x\in A^\mathbb Z : x{[-n,n]}\in L(\sigma)\mbox{ for any }n\}
$$
is a closed $T$-invariant subset of $A^\mathbb Z$, and $T_\sigma$ is the shift $T$ restricted to the set $X_\sigma$.  We will denote the $k$th coordinate of $x \in X_\sigma$ by $x[k]$.

The following theorem is one of the main results in~\cite{durand_host_skau}. A part of this theorem was proved in~\cite{forrest}.

\begin{theorem}[\cite{durand_host_skau}]\label{DHSmain}
The family $\mathfrak{B}$ of Bratteli-Vershik systems associated with stationary,
properly ordered Bratteli diagrams is (up to isomorphism) the disjoint union of
the family of substitution minimal systems and the family of stationary odometer
 systems. Furthermore, the correspondence in question is given by an explicit
and algorithmic effective construction. The same is true of the computation of the
(stationary) dimension group associated with a substitution minimal system.
\end{theorem}

In the proof of Durand, Host and Skau, the explicit construction of homeomorphism is given, while the proof of Forrest has more existential nature. For instance, let $B$ be a stationary properly ordered Bratteli diagram with a simple hat
(a Bratteli diagram has a simple hat whenever it has only simple edges between the top vertex and the first level).
The substitution $\sigma$ \textit{read on} $B$ is defined as follows. Since $B$ is
stationary, all information about it is given by the first level. Let $\mathcal{A} =
\{a_1,\ldots,a_K\}$ be an alphabet, where $K = |V_1| = |V_n|$. With a vertex
number $i$ we associate a letter $a_i$, this operation does not depend on the level.
Consider a letter (vertex) $a \in V_2$ and the ordered list $(e_1, \ldots, e_r)$ of
the edges in $E_1$ with $r(e_j) = a$. Let $(a_1,\ldots, a_r)$ be the ordered list of
 the sources of these edges. Then define $\sigma(a) = a_1 \ldots a_r$. If $\sigma$
is periodic then $(X_B, \varphi_B)$ is isomorphic to an odometer with stationary
base, otherwise $(X_B, \varphi_B)$ is isomorphic to $(X_\sigma, T_\sigma)$.
On the other hand, any substitution $\sigma$ defines a stationary Bratteli diagram: one has to find a stationary ordered Bratteli diagram $B$ with simple hat such that the substitution read on $B$ is exactly $\sigma$.
A substitution $\sigma$ on an alphabet $\mathcal{A}$ is called \textit{proper} if there exists an integer $n > 0$ and two letters $a,b \in \mathcal{A}$ such that for every $c \in \mathcal{A}$, $a$ is the first letter and $b$ is the last letter of $\sigma^n(c)$.
In the case when $\sigma$ is primitive, aperiodic and proper, the systems $(X_\sigma, T_\sigma)$ and $(X_B, \varphi_B)$ are isomorphic. In the case of a non-proper substitution one has to build Bratteli diagram using return words~\cite{durand_host_skau}.

The following theorem extends the result from the minimal case to the aperiodic one.

\begin{theorem}[\cite{BKM09}]\label{TheoremVershikMapIsSubstitution}
(i) Suppose that $(X_B,\varphi_B)$ is an aperiodic Bratteli-Vershik system with $B$ a stationary ordered Bratteli diagram and $X_B$ is perfect. Then the system $(X_B,\varphi_B)$ is  conjugate to an aperiodic substitution dynamical system (with substitution read on $B$) if and only if no  restriction of $\varphi_B$ to a minimal component is isomorphic to an odometer.

(ii) Let $\sigma : A\rightarrow A^+$ be an aperiodic substitution such that $|\sigma^n(a)|\rightarrow \infty$ as $n \rightarrow \infty$ for all $a \in A$. Then the substitution dynamical system $(X_\sigma,T_\sigma)$ is conjugate to the Vershik map of a stationary Bratteli diagram.
\end{theorem}

The technique presented in the proof is applicable for a wide class of substitutions
including various Chacon-like substitutions.

A substitution is {\it recognizable} if for each $x\in X_\sigma$ there exist a unique $y\in X_\sigma$ and unique $i\in \{0,\ldots,|\sigma(y[0])|-1\}$ such that $x = T^i_\sigma\sigma(y)$. Every aperiodic primitive substitution is recognizable~\cite{mosse:1992, mosse:1996}. In~\cite{BKM09}, it is shown that an arbitrary aperiodic substitution is recognizable:

\begin{theorem}[\cite{BKM09}]\label{TheoremRecognizability}
Each aperiodic substitution $\sigma : A\rightarrow A^+$  is recognizable.
\end{theorem}

The proof involves the usage of Downarowicz-Maass' techniques~\cite{DM} (see also Section~\ref{finrank}).

\subsection{Invariant measures on stationary diagrams and their supports}

In this subsection, we give an explicit description of all ergodic probability measures on stationary Bratteli diagrams invariant with respect to the tail equivalence relation (or the Vershik map).

Let $B = (V,E)$ be a Bratteli diagram with incidence matrices ${F_n}$.
For $w \in V_{n}$, the set $E(v_{0}, w)$ defines the clopen subset $X_{w}^{(n)} := \{x = (x_{i}) \in X_{B} : r(x_{n}) = w \}$ of $X_B$. Then  $\{X_{w}^{(n)} : w \in V_{n}\}$ is a clopen partition of $X_{B}$. Analogously, the sets $X_{w}^{(n)}(\overline{e}) := \{x = (x_{i}) \in X_{B} : x_{i} = e_{i}, i = 1,...,n\}$ determine a clopen partition of $X_{w}^{(n)}$ where  $\overline{e} = (e_{1}, \ldots ,e_{n}) \in E(v_{0}, w)$, $n \geq 1$. Note that any two paths $x$, $y$ from $X_{B}$ are $\mathcal{R}$-equivalent if and only if there exists $w \in V$ such that $x \in X_{w}^{(n)}(\overline{e})$ and $y \in X_{w}^{(n)}(\overline{e'})$ for some $\overline{e}, \overline{e'} \in E(v_{0}, w)$. Recall that a measure $\mu$ on $X_B$ is called $\mathcal{R}$-\textit{invariant} if for any two paths $\overline{e}$ and $\overline{e'}$ from $E(v_{0}, w)$ and any vertex $w$, one has $\mu(X_{w}^{(n)}(\overline{e})) = \mu(X_{w}^{(n)}(\overline{e'}))$.  Then
$$
\mu(X_{w}^{(n)}(\overline{e})) =  \frac{1}{h_{w}^{(n)}}\mu(X_{w}^{(n)}), \ \  \overline{e} \in E(v_{0}, w).
$$

For $x = (x_1, \ldots, x_N) \in \mathbb{R}^{N}$, we wil write $x \geq 0$ if $x_i \geq 0$ for all $i$. Let $\mathbb{R}^N_+ = \{x \in \mathbb{R}^N : x \geq 0\}$. Let
$$
C_k^{(n)} = F_k^T \ldots F_n^T(\mathbb{R}^{|V_{n+1}|}_+) \mbox{ for } 1 \leq k \leq n.
$$
Clearly, $\mathbb{R}^{|V_{k}|}_+ \supset C_k^{(n)}\supset C_k^{(n+1)}$ for all $n \geq 1$. Let
$$
C_k^{(\infty)} = \bigcap_{n \geq k}C_k^{(n)} \mbox{ for } k \geq 1.
$$

The following theorem describes a construction of invariant measures on Bratteli diagram of general form.

\begin{theorem}[\cite{BKMS_2010}]\label{Theorem_measures_general_case}  Let $B = (V,E)$ be a Bratteli diagram such that the tail equivalence relation
 $\mathcal R$ on $X_B$ is aperiodic.
If  $\mu\in M(\Rk)$,
then the vectors $p^{(n)}=(\mu(X_w^{(n)}(\ov{e})))_{w\in
V_n}$, $\ov{e}\in E(v_0,w)$, satisfy the following conditions for $n\geq 1$:

(i) $p^{(n)} \in  C_n^\infty$,

(ii) $F_n^T p^{(n+1)} = p^{(n)}$.

Conversely, if a sequence of vectors $\{ p^{(n)}\}$ from
$\mathbb R^{|V_n|}_+$ satisfies condition (ii), then there
exists a non-atomic  finite Borel $\mathcal R$-invariant measure $\mu$ on
$X_B$ with $p_w^{(n)}=\mu(X_w^{(n)}(\overline e))$ for all $n\geq 1$
and $w\in V_n$.

The $\Rk$-invariant measure $\mu$ is a probability measure if and only if

(iii) $\sum_{w\in V_n}h_w^{(n)}p_w^{(n)}=1$ for $n=1$,

\noindent in which case this equality holds for all $n\ge 1$.
\end{theorem}

In \cite{BKMS_2010}, all invariant ergodic measures on a stationary Bratteli diagram were  described as follows. It was first shown  that the study of any stationary Bratteli diagram $B = (V,E)$ (with $|V| = N$) can be reduced to the case when the incidence matrix $F$ of size $N\times N$ has the Frobenius Normal Form:

\begin{equation}\label{Frobenius Form}
F =\left(
  \begin{array}{ccccccc}
    F_1 & 0 & \cdots & 0 & 0 & \cdots & 0 \\
    0 & F_2 & \cdots & 0 & 0 & \cdots & 0 \\
    \vdots & \vdots & \ddots & \vdots & \vdots & \cdots& \vdots \\
    0 & 0 & \cdots & F_s & 0 & \cdots & 0 \\
    X_{s+1,1} & X_{s+1,2} & \cdots & X_{s+1,s} & F_{s+1} & \cdots & 0 \\
    \vdots & \vdots & \cdots & \vdots & \vdots & \ddots & \vdots \\
    X_{m,1} & X_{m,2} & \cdots & X_{m,s} & X_{m,s+1} & \cdots & F_m \\
  \end{array}
\right)
\end{equation}

We can telescope the diagram $B$ and regroup the vertices in such a way that the matrix $F$ will have the following property:
\begin{eqnarray} \label{prop-I}
& F & \mbox{has the form (\ref{Frobenius Form}) where every nonzero matrix $F_i$}\\[-1ex]
 &  & \mbox{on the main diagonal  is primitive.} \nonumber
\end{eqnarray}
We can telescope the diagram $B$ further to make sure that
\begin{eqnarray} \label{prop-I'}
& F & \mbox{has the form (\ref{Frobenius Form}) where every nonzero matrix $F_i$}\\[-1ex]
& & \mbox{on the main diagonal  is strictly positive.} \nonumber
\end{eqnarray}

The matrices $F_i$ determine the partition of the vertex set $V$ into subsets (classes) $V_i$ of vertices. In their turn, these subsets generate  subdiagrams $B_i$. The non-zero matrices $X_{j,k}$ indicate which subdiagrams  are linked by some edges (or finite paths). Notice that each subdiagram $B_i$, $i=1,...,s$, corresponds to a minimal component of the cofinal equivalence relation $\mathcal R$.

We denote by $F_\alpha$, $\alpha \in \Lambda$, the \textit{non-zero} matrices on the main diagonal in (\ref{Frobenius Form}).
Let $\alpha \geq \beta$. It is said that the class of vertices $\alpha$ {\it has access} to a class $\beta$, in symbols $\alpha \succeq \beta$, if and only if either $\alpha = \beta$ or there is a finite path in the diagram from a vertex which belongs to $\beta$ to a vertex from  $\alpha$. In other words, the matrix $X_{\alpha, \beta}$ is non-zero. A class $\alpha$ is called {\it final (initial)} if there is no class $\beta$ such that $ \alpha \succ   \beta$ ($\beta \succ \alpha$).

Let  $\rho_\alpha$ be the spectral radius of $F_\alpha$. A class $\alpha\in \{1,...,m\}$ is called  {\it distinguished} if $\rho_\alpha > \rho_\beta$ whenever $ \alpha \succ \beta$. Notice that all classes $\alpha = 1,\ldots,s$ are necessarily distinguished. A real number $\lambda$ is called a {\it
distinguished eigenvalue} if there exists a non-negative eigenvector $x$ with $Fx=\lambda x$.
A real number $\lambda$ is a distinguished eigenvalue if and only if there exists a distinguished class $\alpha$ such that $\rho_\alpha = \lambda$. In this case we denote $\lambda_\alpha = \lambda$.
If $x = (x_1,...,x_N)^T$ is an eigenvector corresponding to a distinguished eigenvalue $\lambda_\alpha$, then $x_i >0$ if and only if $i \in \beta$ and $\alpha \succeq \beta$.

The main result of \cite{BKMS_2010} completely describes the simplex of $\mathcal R$-invariant probability measures
of a stationary Bratteli diagram. Denote
$$
core (A) = \bigcap_{k \geq 1} A^{k}(\mathbb{R}_+^N).
$$
Let $\{\xi_1,\ldots,\xi_k\}$ be the extreme vectors of the cone $core(A)$. Normalize each vector so that $\sum_{w \in V_1} h_w^{(1)} (\xi_i)_w = 1$. Set $D = \{x \in core(A) : \sum_{w \in V_1} h_w^{(1)} x_w = 1 \}$.

\begin{theorem}[\cite{BKMS_2010}]\label{Theorem_MeasuresErgodicMeasures}
 Suppose that $B$ is a stationary Bratteli diagram such that the
tail equivalence relation $\Rk$ is aperiodic and
the incidence matrix $F$ satisfies (\ref{prop-I}).
Then there is a one-to-one correspondence between vectors $p^{(1)}\in D$ and $\mathcal R$-invariant probability
measures on $X_B$. This correspondence is given by the rule $\mu \leftrightarrow
p^{(1)}= (\mu(X_w^{(1)})/h_w^{(1)})_{w\in V_1}$. Furthermore, ergodic measures correspond to the extreme vectors $\{\xi_1,\ldots,\xi_k\}$. In particular, there exist exactly $k$ ergodic measures.
\end{theorem}

More precisely, fix a distinguished eigenvalue $\lambda$ and let  $x =  (x_1,...,x_N)^T$ be the probability non-negative eigenvector  corresponding to  $\lambda$. Then the ergodic probability  measure $\mu$ defined by $\lambda$ and $x$ satisfies the relation:
\begin{equation}\label{mucount}
\mu(X_{i}^{(n)}(\overline{e})) = \frac{x_i}{\lambda^{n - 1}},
\end{equation}
where $i\in V_n$ and $\overline e$ is a finite path with $s(\overline e) =i$. Therefore, the clopen values set for $\mu$ has has the form:
\begin{equation}\label{formulaSmu}
S(\mu) = \left\{\sum_{i = 1}^{N} k^{(n)}_{i}\frac{x_i}{\lambda^{n - 1}} : 0  \leq k^{(n)}_{i} \leq h^{(n)}_{i}; \; n = 1, 2, \ldots \right\}.
\end{equation}

If $\lambda$ is a non-distinguished Perron-Frobenius eigenvalue for $A$, then the corresponding class of vertices is non-distinguished and $\mathcal R$-invariant measure on $X_B$ is infinite. The following result describes  infinite $\Rk$-invariant measures for stationary Bratteli diagrams.

\begin{theorem} [\cite{BKMS_2010}]\label{th-infmeas}
 Suppose that $B$ is a stationary Bratteli diagram such that the
tail equivalence relation $\Rk$ is aperiodic and
the incidence matrix $F$ satisfies (\ref{prop-I}).
Then the set of ergodic infinite ($\sig$-finite) invariant measures, which are positive
and finite on at least one open set (depending on the measure), modulo a
constant multiple, is in 1-to-1 correspondence with the set of non-distinguished eigenvalues of $A=F^T$.
\end{theorem}

The following theorem gives necessary and sufficient conditions under which a measure on $X_B$ will be $\mathcal{R}$-invariant.

\begin{theorem}[\cite{BKMS_2010}]\label{TheoremStationaryMeasures}
 Suppose that $B$ is a stationary Bratteli diagram such that the
tail equivalence relation $\Rk$ is aperiodic and
the incidence matrix $F$ satisfies (\ref{prop-I}).
Let $\mu$ be a finite Borel $\mathcal R$-invariant measure on $X_B$.
Set $p^{(n)}=(\mu(X_w^{(n)}(\overline e)))_{w\in V_n}$ where  $\overline e\in E(v_0,w)$. Then for the  matrix $A = F^T$ associated to $B$ the following properties hold:

(i) $p^{(n)}=A p^{(n+1)}$ for every $n\geq 1$;

(ii) $ p^{(n)}\in core(A),\ n\geq 1$.

Conversely, if a sequence of vectors $\{ p^{(n)}\}$ from
$\mathbb R^{N}_+$ satisfies condition (ii), then there
exists a finite Borel $\mathcal R$-invariant measure $\mu$ on
$X_B$ with $p_w^{(n)}=\mu(X_w^{(n)}(\overline e))$ for all $n\geq 1$
and $w\in V_n$.

The $\Rk$-invariant measure $\mu$ is a probability measure if and only if

(iii) $\sum_{w\in V_n}h_w^{(n)}p_w^{(n)}=1$ for $n=1$,

\noindent in which case this equality holds for all $n\ge 1$.
\end{theorem}

Let $\sigma : \mathcal A\to \mathcal A^+$ be an aperiodic substitution such that $|\sigma^n(a)| \to \infty$ as $n\to \infty$ for any $a\in \mathcal A$. Denote by $(X_\sigma, T_\sigma)$ the corresponding substitution dynamical system. Let $M_\sigma = (m_{ab})$ be the matrix of substitution $\sigma$, i.e. $m_{ab} = L_a(\sigma(b))$ where $L_a(\sigma(b))$ is the number of $a$ occurring in $\sigma(b)$. Denote by $B_\sigma$ the stationary Bratteli diagram ``read on substitution''. This means that $M_\sigma$ is  the matrix transpose to the incidence matrix of $B_\sigma$. It follows from Theorem \ref{TheoremVershikMapIsSubstitution} that there exists a stationary ordered Bratteli diagram $B(X_\sigma, T_\sigma) = B$ whose Vershik map $\varphi_B$ is conjugate to $T_\sigma$. Thus, we have two Bratteli diagrams associated to $(X_\sigma, T_\sigma)$. Note that the diagram $B$ may have considerably more vertices than the diagram $B_\sigma$.

\begin{theorem}[\cite{BKMS_2010}]\label{measures_and_substitutions}
There is a one-to-one correspondence $\Phi$ between the set of ergodic $T_\sigma$-invariant  probability measures on the space
$X_\sigma$ and the set of ergodic $\mathcal R$-invariant probability measures on the path space $X_{B_\sigma}$ of the stationary diagram $B_\sigma$ defined by substitution $\sigma$. The same statement holds for non-atomic infinite invariant measures.
\end{theorem}

\subsection{Good measures on stationary Bratteli diagrams}

Two probability measures $\mu$ and $\nu$ defined on Borel subsets of a topological space $X$ are called \textit{homeomorphic} or \textit{topologically equivalent} if there exists a self-homeomorphism $h$ of  $X$  such that $\mu = \nu\circ h$, i.e. $\mu(E) = \nu(h(E))$ for every Borel subset $E$ of $X$. In such a way, the set of all Borel probability measures on $X$ is partitioned into equivalence classes. The classification of Borel measures with respect to homeomorphisms started in the paper by Oxtoby and Ulam~\cite{Oxt-Ul}, where they found a criterion for a measure to be homeomorphic to Lebesgue measure on the interval $[0,1]$. The classification of measures on Cantor sets started with classification of Bernoulli measures (see \cite{Navarro, Navarro-Oxtoby} and later papers \cite{Austin, D-M-Y, Akin3, Yingst}). Akin began a systematic study of homeomorphic measures on a Cantor space \cite{Akin1, Akin2}. It was  noted in  \cite{Akin1}  that  there exist continuum  classes of equivalent full non-atomic probability measures on a Cantor set. This fact is based on the existence of a countable base of clopen subsets of a Cantor set. Recall that for a measure $\mu$ on a Cantor space $X$, the \textit{clopen values set} $S(\mu) = \{\mu(U):\,U\mbox{ is clopen in } X\}$.
The set $S(\mu)$ is a countable dense subset of the unit interval, and this set provides an invariant for topologically equivalent measures, although it is not a complete invariant, in general. But for the class of the so called good measures, $S(\mu)$ \textit{is} a complete invariant.
\begin{definition}
A full non-atomic probability measure $\mu$ is \textit{good} if whenever $U$, $V$ are clopen sets with $\mu(U) < \mu(V)$, there exists a clopen subset $W$ of $V$ such that $\mu(W) = \mu(U)$.
\end{definition}
It turns out that such measures are exactly invariant measures of uniquely ergodic
minimal homeomorphisms of Cantor sets (see \cite{Akin2},
\cite{glasner_weiss:1995-1}).
It is obvious that all above mentioned definitions (clopen values set, good measures,  etc.) are applicable to the measures on Bratteli diagrams. In particular, we note that the ergodic measures corresponding to minimal components are automatically good: on a simple stationary Bratteli diagram any Vershik map is minimal and uniquely ergodic.

The following results deal with the classification of ergodic probability Borel measures on stationary non-simple Bratteli diagrams which are invariant with respect to the tail equivalence relation $\mathcal{R}$.
Let $B = (V, E)$ be a stationary non-simple Bratteli diagram, $F$ be its incidence $K \times K$ matrix and $A = F^T$. Let $\mu$ be the measure defined by a distinguished class of vertices $\alpha$ and $\lambda$ the corresponding distinguished eigenvalue of $A$. Denote by  $(y_1,...,y_K)^T$  the probability  eigenvector of the matrix $A$ corresponding to $\lambda$. Notice that the vector $(y_1,...,y_K)^T$ may have zero entries. These zero entries are assigned to the vertices from $B$ that are not accessible from the class $\alpha$. Denote by $(x_1, \ldots, x_n)^T$ the positive vector obtained from $(y_1,...,y_K)^T$ by crossing out zero entries.  We call  $(x_1,...,x_n)^T$ the \textit{reduced vector} corresponding to the measure $\mu$.
The following theorem gives a criterion for $\mu$ to be good.

\begin{theorem}[\cite{S.B.O.K.}]\label{KritGood}
Let $\mu$ be an ergodic $\mathcal R$-invariant measure on a stationary diagram $B$ defined by a distinguished eigenvalue $\lambda$ of the matrix $A = F^T$. Denote by $x = (x_1,...,x_n)^T$ the corresponding reduced vector. Let the vertices $m+1, \ldots, n$ belong to the distinguished class $\alpha$ corresponding to $\mu$. Then $\mu$ is good if and only if   there exists $R \in \mathbb{N}$ such that $\lambda^R x_1,...,\lambda^R x_m$ belong to the additive group generated by $\{x_j\}_{j=m+1}^n$.

If the clopen values set of $\mu$ is rational and $(\frac{p_1}{q}, \ldots, \frac{p_n}{q})^T$ is the corresponding reduced vector, then $\mu$ is good if and only if $\gcd(p_{m+1}, ...,p_n) | \; \lambda^R$ for some $R \in \mathbb N$.
\end{theorem}

The idea of the proof is as follows. The support of measure $\mu$ is a stationary simple subdiagram $B_\alpha$ which corresponds to the distinguished class $\alpha$. The measure $\mu\mid_{B_\alpha}$ is good since it is a unique ergodic invariant measure for a Vershik homeomorphism. The property of goodness can be destroyed when $\mu\mid_{B_\alpha}$ is extended to the measure $\mu$ on the whole diagram $B$. The measures of cylinder sets that end in the vertices of $B$ that do not belong to class $\alpha$ are obtained as infinite sums of measures of cylinder sets that end in the vertices of class $\alpha$. Thus, new values in $S(\mu)$ might appear which do not belong to $S(\mu\mid_{B_\alpha})$, and this will be the reason for $\mu$ to be not good.

For a good measure $\mu$, there always exists a subgroup $G \subset \mathbb R$ such that $S(\mu) = G\cap [0,1]$, i.e. $S(\mu)$ is \textit{group-like} (see~\cite{Akin2}). Indeed, let $U$ and $V$ be any two clopen sets and $\mu$ be a good measure. Suppose $\mu(U) = \alpha$, $\mu(V) = \beta$ and without loss of generality $\alpha < \beta$. Thus, $\alpha, \beta \in S(\mu)$. Since $\mu$ is good, there is a clopen set $W \subset V$ such that $\mu(W) = \alpha$. Hence $V \setminus W$ is a clopen set of measure $\beta - \alpha$, and $\beta - \alpha \in S(\mu)$. It easily follows that $S(\mu)$ is group-like.

It turns out that for any ergodic invariant measure $\mu$ on a stationary diagram $B$ the set $S(\mu)$ is group-like~\cite{S.B.O.K.}. Moreover, it was proved that for an ergodic invariant measure $\mu$ defined by a distinguished eigenvalue $\lambda$ and distinguished class of vertices $\alpha$, the following equation holds:
$$
S(\mu) = \left(\bigcup_{N=0}^\infty \frac{1}{\lambda^N} H\right) \cap [0,1],
$$
where $(x_1, \ldots, x_n)^T$ is the corresponding reduced vector and $H$ is the additive subgroup of $\mathbb{R}$ generated by $\{x_1, \ldots , x_n\}$.
Two cases were considered in the proof: (i) $\lambda$ is rational and hence $S(\mu)\subset\mathbb Q$; (ii) $\lambda$ is an irrational algebraic integer of degree $k$, then $S(\mu)\subset \mathbb{Q}[\lambda]$, where $\mathbb{Q}[\lambda]$ is the least ring that contains both $\mathbb{Q}$ and $\lambda$. The case (i) is relatively simple. The following asymptotics mentioned in \cite{BKMS_2010} was used:
$$
(A^n)_{ij} \sim \lambda^n,\ \ n\to \infty,\ \ \mbox{for}\ i\in \beta,\ j\in \alpha,\ \mbox{with}\
\alpha \succeq \beta.
$$
Here $\sim$ means that the ratio tends to a positive constant. On the other hand,
$$
(A^n)_{ij} = o(\lambda^n),\ \ n\to \infty,\ \ \mbox{for}\ j\in
\beta \prec \alpha.
$$
In the case (ii), the eigenvector entries and eigenvalues are represented as vectors with rational entries (the dimension of vectors is equal to the algebraic degree of $\lambda$). Then the operation of dividing by $\lambda$ is a linear transformation in the vector space $\mathbb{Q}^{k}$. The methods of linear algebra and matrix theory were used to obtain the proof.

Theorem~\ref{KritGood} was used to prove the following result. We showed that one can build infinitely many homeomorphic ergodic invariant measures on stationary diagrams such that the corresponding tail equivalence relations are non-orbit equivalent.
\begin{theorem}[\cite{S.B.O.K.}]\label{goodhomeo}
Let $\mu$ be a good ergodic $\mathcal{R}$-invariant probability measure on a stationary (non-simple) Bratteli diagram $B$. Then there exist stationary Bratteli diagrams  $\{B_i\}_{i=0}^\infty$ and  good ergodic $\mathcal{R}_i$-invariant probability measures $\mu_i$ on $B_i$ such that each measure $\mu_i$ is homeomorphic to $\mu$ and the dynamical systems $(B_i, \mathcal R_i)$, $(B_j, \mathcal R_j)$ are topologically orbit equivalent if and only if $i = j$. Moreover, the diagram $B_i$ has exactly $i$ minimal components for the tail equivalence relation $\mathcal R_i, i\in \mathbb N$.
\end{theorem}

In~\cite{K12}, the notion of good measure is extended to the case of infinite measures on Cantor sets.
Borel infinite measures arise as ergodic invariant measures for aperiodic
homeomorphisms of a Cantor set. The study of homeomorphic infinite measures is
of crucial importance for the classification of Cantor aperiodic systems up to orbit equivalence.

The direct analogues of Theorems~\ref{KritGood}, \ref{goodhomeo} for infinite measures
on stationary Bratteli diagrams can be found in~\cite{K12}. Measures on stationary
Bratteli diagrams can be also considered as extensions of good measures on
 non-compact locally compact Cantor sets (see~\cite{K122} for details).

In the paper \cite{BH}  the authors initiated the study of  properties of traces of dimension groups which are motivated by Cantor dynamics. In particular, it was shown that there are example of minimal Cantor systems whose ergodic measures have completely different (in some sense even opposite) properties.

\subsection{Complexity and orbit equivalence}
Let $\mu$ be an ergodic invariant measure for a proper primitive substitution dynamical system. By Theorem~\ref{DHSmain}, we can find easily the clopen values set $S(\mu)$ in terms of the matrix of substitution. By Corollary~\ref{OE corollary}, to construct a minimal substitution dynamical system which is orbit equivalent to a given one, it suffices to find another stationary simple Bratteli diagram such that the clopen values set is kept unchanged. In the paper~\cite{BK14}, two constructions are used to build countably many non-isomorphic orbit equivalent minimal substitution dynamical systems. In both constructions, the complexity function $n\mapsto p_\sigma(n)$ is used to distinguish non-isomorphic systems. Recall that the function $p_\sigma(n)$ counts the number of words of length $n$ in the infinite sequence invariant with respect to $\sigma$. If minimal substitution dynamical systems $(X_\sigma, T_\sigma)$ and $(X_\zeta, T_\zeta)$ are topologically conjugate, then there exists a constant $c$ such that, for all $n > c$ one has $p_\sigma(n - c) \leq p_\zeta(n) \leq p_\sigma(n + c)$ (see~\cite{F}). In the proof of the following theorem, the complexity of the systems $(X_{\zeta_n}, T_{\zeta_n})$ is forced to grow by increasing the number of letters in the alphabets of $\zeta_n$.

\begin{theorem}[\cite{BK14}]\label{diff_alphabets}
Let $\sigma$ be a proper substitution.
Then there exist countably many proper substitutions $\{\zeta_n\}_{n=1}^\infty$ such that $(X_\sigma, T_\sigma)$ is orbit equivalent to $(X_{\zeta_n}, T_{\zeta_n})$, but the systems $\{(X_{\zeta_n}, T_{\zeta_n})\}_{n=1}^\infty$ are pairwise non-isomorphic.
\end{theorem}

Moreover, if one wants to have a substitution dynamical system which is strongly orbit equivalent to a given $(X_\sigma, T_\sigma)$, then additionally the dimension group of the diagram $B_\sigma$ must be unchanged. If $A$ is the incidence matrix for $\sigma$ and $A^N$ is the incidence matrix for $\zeta$ for some $N \in \mathbb{N}$, then the dimension groups associated to minimal Cantor systems $(X_\sigma, T_\sigma)$ and $(X_\zeta, T_\zeta)$ are order isomorphic (see~\cite{giordano_putnam_skau:1995}). In the following theorem, incidence matrices of built substitution systems are the powers of the incidence matrix of the initial substitution system $(X_\sigma, T_\sigma)$. Hence substitutions $\sigma$ and $\zeta_n$ have the same fixed alphabet, and the complexity function is made increasing by enlarging the length of substitution $\zeta_n$ and by an appropriate permutation of letters. Thus, we produce a countable family of pairwise non-isomorphic strong orbit equivalent substitution systems.

\begin{theorem}[\cite{BK14}]
Let $\sigma$ be a primitive proper substitution. Let $(B, \leq)$ be the corresponding stationary properly ordered simple Bratteli diagram.
Then there exist countably many telescopings $B_n$ of $B$ with proper orders $\leq_n$ and corresponding substitutions $\zeta_n$ read on $B_n$ such that the substitution dynamical
systems $\{(X_{\zeta_n}, T_{\zeta_n})\}_{n=1}^\infty$ are pairwise non-isomorphic and strong orbit equivalent to $(X_\sigma, T_\sigma)$.
\end{theorem}

In the following theorem, given a Bratteli-Vershik system $(X_B, \varphi_B)$ on a stationary simple diagram, an orbit equivalent Bratteli-Vershik system is found with the least possible number of vertices. The vector technique developed in~\cite{S.B.O.K.} is used to prove the result. Let $\lambda$ be a Perron-Frobenius eigenvalue for the incidence matrix of $B$. Since the algebraic degree $\deg \lambda$ is equal to $k$, the dimension of the vector space corresponding to any other orbit equivalent Bratteli-Vershik system is at least $k$. Hence there is no stationary Bratteli-Vershik system with less than $k$ vertices which is orbit equivalent to $(X_B, \varphi_B)$.

\begin{theorem}[\cite{BK14}]\label{minBD}
Let $\sigma$ be a primitive substitution whose incidence matrix has a Perron-Frobenius eigenvalue $\lambda$ and $k = \deg \lambda$. Then $(X_\sigma, T_\sigma)$ is orbit equivalent to a Bratteli-Vershik system defined on a stationary Bratteli diagram with $k$ vertices on each level. Moreover, there is no stationary Bratteli-Vershik system with less than $k$ vertices which is orbit equivalent to $(X_\sigma, T_\sigma)$.
\end{theorem}

\section{Finite rank Bratteli diagrams}\label{finrank}
This section is devoted to the study of aperiodic Cantor dynamical systems which can be represented by Bratteli diagrams with uniformly bounded number of vertices on each level. It is an open question which classes of Cantor dynamical systems admit such a representation.

\begin{definition}
A Cantor dynamical system $(X,S)$ has the {\it topological rank} $K>0$ if it admits a Bratteli-Vershik model $(X_B,\varphi_B)$ such that the number of vertices of the diagram $B$ at each level is not greater than $K$ and $K$ is the least possible number of vertices for any Bratteli-Vershik realization.
\end{definition}

If a system $(X,S)$ has the rank $K$, then, by an appropriate telescoping, we can assume that the diagram $B$ has exactly $K$ vertices at each level. It is said that a homeomorphism $S: X\rightarrow X$ is {\it expansive} if there exists $\delta>0$ such that for any distinct $x,y\in X$ there is $m\in\mathbb Z$ with $d(S^mx,S^my)>\delta$. The number $\delta$ is called an {\it expansive constant}. Note that the notion of expansiveness does not depend on the choice of the metric $d$, see \cite{walters:book}. The series of papers is devoted to the characterization of Cantor minimal systems through expansive Bratteli diagrams, namely, expansive diagrams with constant number of incoming edges characterize Toeplitz systems~\cite{gjerde_johansen:2000}, expansive diagrams with a finite set of incidence matrices characterize linearly recurrent subshifts~\cite{CDHM03}. The main result in~\cite{DM} states that

\begin{theorem}[\cite{DM}]\label{minfinranCSexpansive}
Every Cantor minimal system of finite rank $d > 1$ is expansive.
\end{theorem}

In~\cite{DM}, the authors suggested a method of coding the 
dynamics on Bratteli diagrams by means of the so-called  $j$-symbols. 
A new proof of Theorem \ref{minfinranCSexpansive} is given in \cite{H14}.
In~\cite{BKM09}, the ideas and results from~\cite{DM} are 
generalized, and the abstract definitions of $j$-symbols, $j$-sequences
 etc. are given. They can be used to study dynamics of different nature, 
for instance, Bratteli-Vershik systems and substitution dynamical 
systems. The main advantage of this approach is that it allows one to 
use the machinery of symbolic dynamics for solving some problems of 
Cantor dynamics.
The following result is an extension 
 of Theorem~\ref{minfinranCSexpansive} and was proved by using the 
technique of $j$-sequences.

\begin{theorem}[\cite{BKM09}]\label{TheoremFiniteRank}
Let $(X,S)$  be an aperiodic Cantor dynamical system of finite rank. If the restriction of $(X,S)$ to every minimal component is not  conjugate to an odometer, then $(X,S)$ is expansive.
\end{theorem}

In~\cite{BKMS_2013}, the structure of invariant measures on finite rank Bratteli diagrams is considered. In particular, it is shown that every ergodic invariant measure (finite or ``regular'' infinite) can be obtained as an extension from a simple vertex subdiagram.

\begin{definition}
Let $B$ be a Bratteli diagram. By a {\it vertex subdiagram} of $B$, we mean a Bratteli diagram $\ov B = (\ov V,\ov E)$ constructed by taking some vertices at each level $n$ of the diagram $B$ and then considering {\it all} the edges of $B$ that connect these vertices.
\end{definition}

Let $\ov B = (\ov V,\ov E)$ be a subdiagram of $B$. Consider the set $X_{\ov B}$ of all infinite paths of the subdiagram $\ov B$. Then the set $X_{\ov B}$ is naturally seen as a subset of $X_B$. Let $\ov \mu$ be a finite invariant (with respect to the tail equivalence relation $\mathcal R$) measure on $X_{\ov B}$. Let $\wh X_{\ov B}$ be the saturation of $X_{\ov B}$ with respect to $\mathcal R$. In other words, a path $x\in X_B$ belongs to $\wh X_{\ov B}$ if it is $\mathcal R$-equivalent to a path $y\in X_{\ov B}$. Then $\wh X_{\ov B}$ is $\mathcal R$-invariant and $X_B$ is a complete section for $\mathcal R$ on $\wh X_{\ov B}$. By the {\it extension of measure $\ov \mu$ to $\wh X_{\ov B}$} we mean the $\mathcal R$-invariant measure $\widehat {\ov \mu}$ on $\wh X_{\ov B}$ (finite or infinite) such that $\widehat {\ov\mu}$ induced on $X_{\ov B}$ coincides with $\ov \mu$. To extend the measure $\widehat {\ov\mu}$ to the $\mathcal R$-invariant measure on the whole space $X_B$, we set $\widehat {\ov\mu}( X_B \setminus \wh X_{\ov B}) = 0$.

In the next theorem, we describe the structure of the supports of ergodic invariant measures. Everywhere below the term ``measure'' stands for an $\mathcal R$-invariant measure. By an infinite measure we mean any $\sigma$-finite non-atomic measure which is finite (non-zero) on some clopen set. The support of each ergodic measure turns out to be the set of all paths that stabilize in some subdiagram, which geometrically can be seen as ``vertical'', i.e., they will eventually stay in the subdiagram. Furthermore, these subdiagrams are pairwise disjoint for different ergodic measures.  It is shown in~\cite{BKMS_2013}, that for any finite rank diagram $B$ one can find finitely many vertex subdiagrams $B_\alpha$ such that each finite ergodic measure on $X_{B_\alpha}$ extends to a (finite or infinite) ergodic measure on $X_B$. It is also proved that each ergodic measure (both finite and infinite) on $X_B$ is obtained as an extension of a finite ergodic measure from some $X_{B_\alpha}$. Moreover, the following theorem holds:

\begin{theorem}[\cite{BKMS_2013}]\label{TheoremGeneralStructureOfMeasures}
Let $B$ be a Bratteli diagram of finite rank $d$.
The diagram $B$ can be telescoped in such a way that for every
probability ergodic measure $\mu$ there exists a subset $W_\mu$ of
vertices from $\{1,\ldots,d\}$ such that the support of $\mu$ consists of
all infinite paths that eventually go along the vertices of $W_\mu$
only. Furthermore,

(i) $W_\mu\cap W_\nu = \emptyset$ for different ergodic measures
$\mu$ and $\nu$;

(ii) given a probability ergodic measure $\mu$, there exists a constant $\delta>0$
such that for any $v\in W_\mu$ and any level $n$
$$\mu(X_v^{(n)})\geq \delta
$$
where $X_v^{(n)}$ is the set of all paths that go through the vertex $v$ at level $n$;

(iii) the subdiagram generated by $W_\mu$ is simple and uniquely
ergodic. The only ergodic measure on the path space of the subdiagram
is the restriction of measure $\mu$.

If a probability ergodic measure $\mu$ is the extension of a measure
 from the vertical subdiagram determined by a proper subset $W\subset\{1,\ldots,d\}$, then
$$\lim_{n\to\infty}\mu(X_v^{(n)}) = 0\mbox{ for all }v\notin W.$$
\end{theorem}

We observe that condition (ii) can be used in practice to determine the support 
of  an ergodic measure $\mu$.

A number of sufficient conditions for unique ergodicity of a finite rank Bratteli-Vershik system are obtained in~\cite{BKMS_2013}. Here we present one of them.

\begin{definition}
(i) For two positive vectors $x,y\in\mathbb R^d$, the {\it projective metric (Hilbert metric)} is
$$
D(x,y) = \ln\max_{i,j}\frac{x_iy_j}{x_jy_i} = \ln \frac{\max_i\frac{x_i}{y_i}}{\min_j \frac{x_j}{y_j}},
$$
where $(x_i)$ and $(y_i)$ are entries of the vectors $x$ and $y$.

(ii) For a non-negative matrix $A$, the {\it Birkhoff contraction
coefficient} is $$\tau(A) = \sup_{x,y>0}\frac{D(Ax,Ay)}{D(x,y)}.$$
\end{definition}

\begin{theorem}[\cite{BKMS_2013}]\label{theoremUniqueErgodicityInTermsTau}
Let $B$ be a simple Bratteli diagram of finite rank
with incidence matrices $\{F_n\}_{n\geq 1}$. Let $A_n = F_n^T$.
Then the diagram $B$ is uniquely ergodic if and only if 
$$\lim_{n\to
\infty}\tau(A_m\ldots A_n) = 0\mbox{ for every }m.
$$
\end{theorem}

For a positive matrix $A = (a_{i,j})$, set
$$\phi(A) = \min_{i,j,r,s} \frac{a_{i,j}a_{r,s}}{a_{r,j}a_{i,s}}.
$$
If $A$ has a zero entry, then, by definition, we put $\phi(A)=0$.
As noticed in \cite{hartfiel:book:2002}, 
 $$\tau(A) = \frac{1-\sqrt{\phi(A)}}{1+\sqrt{\phi(A)}}$$
when  $A$ has a nonzero entry in each row.

The following results give computable sufficient conditions for measure 
uniqueness.

\begin{proposition}\label{PropositionSufficientConditionsUniqueErgodicity}
Let $\{A_n\}_{n\geq 1} = F^T_n$ be primitive incidence matrices of a finite rank diagram $B$.

(1) If $$\sum_{n = 1}^\infty \sqrt{\phi(A_n)} = \infty,$$
then $B$ admits a unique invariant probability measure.

(2) If
$$
\sum_{n = 1}^\infty \left(\frac{m_n}{M_n} \right)=\infty,
$$
where $m_n$ and $M_n$ are the smallest and the largest entry of
$A_n$ respectively, then $B$ admits a unique invariant probability
measure.

(3)   If  $||F_n||_1\leq Cn$ for some $C>0$ and all sufficiently
large $n$, then the diagram admits a unique invariant probability
measure\footnote{Here $||A||_1 = \sum_{i,j}|a_{i,j}|$.}. In particular, this result holds if the diagram has only finitely
many different incidence matrices.

\end{proposition}

\begin{example} Let $B$ be a simple Bratteli diagram with incidence matrices
$$
F_n = \left(
  \begin{array}{cccc}
    f_1^{(n)} & 1 & \cdots & 1 \\
    1 & f_2^{(n)} & \cdots & 1 \\
    \vdots & \vdots & \ddots & \vdots \\
    1 & 1 & \cdots & f_d^{(n)} \\
      \end{array}
\right).
$$
Let $q_n = \mbox{max}\{f_i^{(n)}f_j^{(n)} : i\neq j \}$. If for $A_n = F_n^T$
$$
\sum_{n=1}^\infty \sqrt{\phi(A_n)} = \sum_{n=1}^\infty
\frac{1}{\sqrt q_n} =\infty,
$$
then there is a unique invariant probability measure on $B$. 
\end{example}

We say that a Bratteli diagram of a finite rank is of {\it exact finite rank} if there is a finite invariant measure $\mu$ and a constant $\delta >0$ such that after a telescoping $\mu(X_v^{(n)})\geq \delta$ for all levels $n$ and vertices $v$.
The following result shows that the Vershik map on the path space of an exact finite rank diagram cannot be strongly mixing independently of the ordering.

\begin{theorem}[\cite{BKMS_2013}]\label{TheoremAbsenceStrongMixing} Let $B = (V,E,\omega)$
be an ordered simple Bratteli diagram of exact finite
rank. Let $\varphi_{\omega}:X_B\rightarrow X_B$ be the Vershik map defined by the
order $\omega$ on $B$ ($\varphi_{\omega}$ is not necessarily continuous everywhere).
 Then the dynamical system $(X_B,\mu,\varphi_{\omega})$ is not strongly mixing
with respect to the unique invariant measure $\mu$.
\end{theorem}

On the other hand, it is proved in the same paper that for the so-called ``consecutive'' ordering, the Vershik map is not strongly mixing on all finite rank diagrams.

\medskip
The last part of this section is devoted to the following problem. Let $B$ be a  Bratteli diagram of finite rank $k$. It is known (see, e.g., \ref{TheoremGeneralStructureOfMeasures} that $B$ can support at most $k$ ergodic (finite and infinite) measures. Is it possible to determine under what conditions on the incidence matrices of $B$ there exist exactly $k$ ergodic measures? We give a criterion for the existence of $k$ measures when the incidence matrices satisfy the equal row sum property. 

\begin{theorem}[\cite{ABKK}]
Let $B = (V,E)$ be a Bratteli diagram of rank $k \geq 2$; identify $V_n$ with $\{1, ... ,k\}$ for any $n \geq 1$. Let $F_n = (f_{i,j}^{(n)})$ form a sequence of  incidence matrices of $B$ such that
$\sum_{j \in V_n} f_{i,j}^{(n)} = r_n \geq 2$ for every $i \in V_{n+1}$. Suppose that  $\rank\ F_n = k$ for all $n$. Denote
$$
z^{(n)} =
\det
\begin{pmatrix}
\dfrac{f_{1,1}^{(n)}}{r_n} & \ldots & \dfrac{f_{1,{k-1}}^{(n)}}{r_n} & 1\\
\vdots & \ddots & \vdots & \vdots\\
\dfrac{f_{{k},1}^{(n)}}{r_n} & \ldots & \dfrac{f_{{k},{k-1}}^{(n)}}{r_n} & 1
\end{pmatrix}.
$$
Then there exist exactly $k$ ergodic invariant measures on $B$ if and only if
$$
\prod_{n = 1}^{\infty} |z^{(n)}| > 0,
$$
or, equivalently,
$$
\sum_{n = 1}^{\infty}(1 - |z^{(n)}|) < \infty.
$$
\end{theorem}

In the paper \cite{BDM10} the authors addressed the following question:  how can combinatorial properties of Bratteli diagrams be used to describe continuous and measurable eigenvalues?   They used the sequence of incidence matrices to produce necessary (algebraic) conditions for a number $\lambda$ to be an eigenvalue that is both continuous and measurable. It is remarkable that the necessary conditions do not depend on an order of the diagram (that specifies a Vershik map). These results were applied to Toeplitz minimal systems.

\section{Measures and subdiagrams}\label{sbd}

In this section, we consider ergodic $\mathcal R$-invariant measures on arbitrary Bratteli diagrams related to subdiagrams. Suppose that we have a Bratteli diagram $B$ and an ergodic measure $\mu$. It is still an open question whether one can explicitly describe the support of $\mu$ on $X_B$ in terms of the diagram $B$. It would be nice to have a statement similar to \ref{TheoremGeneralStructureOfMeasures}. We give answers to the following questions:

(A) Given a subdiagram $B'$ of $B$ and an ergodic measure $\mu$ on $X_B$, under what conditions on $B'$ the subset $X_{B'}$ has positive measure $\mu$ in $X_B$?

(B) Let $\nu$ be a measure supported by the path space $X_{B'}$ of a  subdiagram $B' \subset B$. Then $\nu$ is extended to the subset $\mathcal R(X_{B'})$ by invariance with respect to the tail equivalence relation $\mathcal R$. Under what conditions $\nu(\mathcal R(X_{B'}))$ is finite (or infinite)?

By a Bratteli subdiagram, we mean a Bratteli diagram $B'$ that can be obtained from $B$ by removing some vertices and edges from each level of $B$. Then $X_{B'}\subset X_B$. We will consider two extreme cases of Bratteli subdiagrams: vertex subdiagram (when we fix a subset of vertices at each level and take all edges between them) and edge subdiagram (some edges are removed from the initial Bratteli diagram but the vertices are not changed). It is clear that an arbitrary subdiagram can be obtained as a combination of these cases.

We keep the following notation:
$\overline{X}_v^{(n)}$ stands for the tower in a subdiagram $\overline{B}$ that is determined by  a vertex $v$ of $\ol B$. Thus,  we consider the paths in $\overline{X}_v^{(n)}$ that contain edges from $\overline{B}$ only. Let $\overline{h}_v^{(n)}$ be the  height of the tower $\ol X_v^{(n)}$. As a rule, objects related to a subdiagram $\ol B$ are denoted by barred symbols.
The following theorem gives criteria for finiteness of the measure extension.

\begin{theorem}[\cite{BKK_14}]\label{thm from BKK}
Let $B$ be a Bratteli diagram with the sequence of incidence matrices $\{F_n\}_{n=0}^\infty$, and let $\ov B$ be a vertex subdiagram of $B$ defined by the sequence of subsets  $\{W_n\}_{n = 0}^\infty, \ W_n \subset V_n$. Suppose that $\ov \mu$ is  a probability  $\ol{\mathcal R}$-invariant measure on $X_{\ov B}$. Then the following properties are equivalent:
\begin{eqnarray*}
\widehat{\overline{\mu}}(\wh X_{\ov B}) < \infty & \Longleftrightarrow & \sum_{n=1}^\infty \sum_{v \in W_{n+1}} \sum_{w \in W'_n} f_{v,w}^{(n)} h_w^{(n)} \ov p_v^{(n+1)} < \infty\\
& \Longleftrightarrow & \sum_{n=1}^{\infty} \sum_{w \in W_{n+1}} \wh{\ov\mu}(X_{w}^{(n+1)})\sum_{v \in W^{'}_{n}} q_{w,v}^{(n)}
< \infty\\
& \Longleftrightarrow & \sum_{i=1}^{\infty} \left(\sum_{w\in W_{i+1}} h_w^{(i+1)} \ov p_w^{(i+1)} - \sum_{w\in W_{i}} h_w^{(i)} \ov p_w^{(i)}\right)< \infty.
\end{eqnarray*}
\end{theorem}

If $(F_n)$ is a  sequence of incidence matrices  of a Bratteli diagram $B$, then we can also define the sequence of stochastic matrices $(Q_n)$ with entries
$$
q^{(n)}_{v,w} = f^{(n)}_{v,w}\frac{h_w^{(n)}}{h_v^{(n+1)}},\ \ v\in V_{n+1},\ w\in V_{n}.
$$
Paper \cite{BKK_14} contains also some necessary and sufficient conditions for finiteness of the measure extension. For instance, it is shown that if
\begin{equation*}\label{suffic}
\sum_{n=1}^{\infty} \sum_{v \in W_{n+1}} \sum_{w \in W^{'}_{n}} q_{v,w}^{(n)} < \infty,
\end{equation*}
then any probability measure $\ov \mu$ defined on the path space $X_{\ol B}$ of the vertex subdiagram $\ol B$ extends to a finite measure $\widehat{\ov \mu}$ on $\wh X_{\ol B}$.

The analogue of Theorem~\ref{thm from BKK} can be proved also for edge subdiagrams:

\begin{theorem}[\cite{ABKK}]\label{edge_fin_krit}
Let $\overline{B}$ be an edge subdiagram of a Bratteli diagram $B$. For a probability invariant measure $\overline{\mu}$ on $X_{\overline{B}}$, the extension $\widehat{\overline{\mu}} (\widehat{X}_{\overline{B}})$ is finite if and only if
$$
\sum_{n=1}^{\infty} \sum_{v \in V_{n+1}} \sum_{w \in V_{n}} \widetilde{f}_{v,w}^{(n)} h_w^{(n)}\overline{p}_v^{(n+1)} < \infty
$$
where $\widetilde{f}_{v,w} = {f}_{v,w} - \ol{f}_{v,w}$.
\end{theorem}

The following theorem gives a necessary and sufficient condition for a  subdiagram $\ol B$ of $B$ to have a path space of zero measure in $X_B$.   Though the theorem is formulated for a vertex subdiagram, the statement remains  true also for any edge subdiagram $\ol B$.

\begin{theorem}[\cite{ABKK}]\label{main}
Let $B$ be a simple Bratteli diagram, and let $\mu$ be any probability ergodic measure on $X_B$. Suppose that  $\ov B$ is a vertex subdiagram of $B$ defined by a sequence  $(W_n)$ of subsets of $V_n$. Then $\mu(X_{\ov B}) = 0$ if and only if
\begin{equation}\label{thin set condition}
\forall \varepsilon > 0 \ \exists n = n(\varepsilon) \ \mathrm{such\ that} \ \forall w \in W_n \ \mathrm{one\ has} \ \frac{\ov h_w^{(n)}}{h_w^{(n)}} < \varepsilon.
\end{equation}
\end{theorem}

In fact, Theorem \ref{main} states that if  a subdiagram $\ol B$ satisfies  (\ref{thin set condition}), then $X_{\ol B}$ has measure zero with respect to every ergodic invariant measure, that is the set $X_{\ol B}$ is {\em thin} according to the definition from \cite{giordano_putnam_skau:2004}. The following result is a corollary of Theorem \ref{main}:

\begin{theorem}[\cite{ABKK}]\label{corol_infty}
Let $\ol B$ be a subdiagram of $B$ such that $X_{\ol B}$ is a thin subset of $X_B$.  Then for any probability invariant measure $\ov \mu$ on $\ov B$ we have $\widehat{\overline{\mu}}(\wh X_{\ov B}) = \infty$.
\end{theorem}

\begin{remark}
There are a lot of papers, where invariant measures for various Bratteli diagrams are studied. For instance, in~\cite{Frick_Petersen, PeVa} the authors consider ergodic invariant probability measures on a Bratteli diagram of a special form, called an Euler graph; the combinatorial properties of the Euler graph are connected to those of Eulerian numbers. The authors of~\cite{FrO} study spaces of invariant measures for a class of dynamical systems which is called polynomial odometers. These are adic maps on regularly structured Bratteli diagrams and include the Pascal and Stirling adic maps as examples. K. Petersen~\cite{Petersen} considers ergodic invariant measures on a Bratteli-Vershik dynamical system, which is based on a diagram whose path counts from the root are the Delannoy numbers.

We would like to mention also the interesting paper by Fisher \cite{F09} where various properties of Bratteli diagrams and measures are discussed. 
\end{remark}

\medbreak
{\bf Acknowledgement.} We are thankful to our colleagues and co-authors for numerous useful discussions of the concept of Bratteli diagrams.  S.B. would like to thank the University of Iowa for the warm hospitality where the paper was finished.

\bibliographystyle{amsalpha}

\end{document}